%            This document is written in Plain TeX

\magnification=1100
\overfullrule0pt

\input amssym.def
\input prepictex
\input pictex
\input postpictex

% ********************* Definitions ************************************

% *****Matrices in subscripts********

\def\QQ{{\Bbb Q}}
\def\RR{{\Bbb R}}
\def\ZZ{{\Bbb Z}}

\def\cT{{\cal T}}

\def\fh{{\frak h}}

\newsymbol\ltimes 226E % This provides the semidirect product symbol
\newsymbol\rtimes 226F % This provides the other semidirect product symbol
%\def\widetilde{\mathaccent"0365 }

% ******************  QED SIGNS  *********************************

\def\qed{\hbox{\hskip 1pt\vrule width4pt height 6pt depth1.5pt \hskip 1pt}}

\def\sqr#1#2{{\vcenter{\vbox{\hrule height.#2pt
\hbox{\vrule width.#2pt height#1pt \kern#1pt
\vrule width.2pt}
\hrule height.2pt}}}}

  % open square

% ********************* FONTS ************************************

\font\smallcaps=cmcsc10
\font\titlefont=cmr10 scaled \magstep2

\font\sectionfont=cmbx10
\font\tinyrm=cmr10 at 8pt

% ******************** SECTION HEADERS ***************************

\newcount\sectno
\newcount\subsectno
\newcount\resultno

\def\section #1. #2\par{
\sectno=#1
\resultno=0
\bigskip%\noindent
\centerline{\sectionfont #1.  #2} }%~\medbreak}

\def\subsection #1\par{\bigskip\noindent{\it  #1} \medbreak}

%******************* MATHEMATICAL LABELS **************************

\def\conj{ \global\advance\resultno by 1
\bigskip\noindent{\bf Conjecture \the\sectno.\the\resultno. }\sl}
\def\prop{ \global\advance\resultno by 1
\bigskip\noindent{\bf Proposition \the\sectno.\the\resultno. }\sl}
\def\lemma{ \global\advance\resultno by 1
\bigskip\noindent{\bf Lemma \the\sectno.\the\resultno. }
\sl}
\def\remark{ \global\advance\resultno by 1
\bigskip\noindent{\bf Remark \the\sectno.\the\resultno. }}
\def\example{ \global\advance\resultno by 1
\bigskip\noindent{\bf Example \the\sectno.\the\resultno. }}
\def\cor{ \global\advance\resultno by 1
\bigskip\noindent{\bf Corollary \the\sectno.\the\resultno. }\sl}
\def\thm{ \global\advance\resultno by 1
\bigskip\noindent{\bf Theorem \the\sectno.\the\resultno. }\sl}
\def\defn{ \global\advance\resultno by 1
\bigskip\noindent{\it Definition \the\sectno.\the\resultno. }\slrm}
\def\endthm{\rm\bigskip}

\def\pf{\rm\bigskip\noindent{\it Proof. }}
\def\endpf{\qed\hfil\bigskip}

%*************** EQUATIONS WITH NUMBERS **************

\def\formula{\global\advance\resultno by 1
\eqno{(\the\sectno.\the\resultno)}}
\def\formulano{\global\advance\resultno by 1 (\the\sectno.\the\resultno)}
\def\tableno{\global\advance\resultno by 1
\the\sectno.\the\resultno. }
\def\lformula{\global\advance\resultno by 1
\leqno(\the\sectno.\the\resultno)}

%************Commutative diagrams**********************

\def\mapright#1{\smash{\mathop
        {\longrightarrow}\limits^{#1}}}

\def\mapdown#1{\Big\downarrow
   \rlap{$\vcenter{\hbox{$\scriptstyle#1$}}$}}

%********** DATING ******************************************
\def\monthname {\ifcase\month\or January\or February\or March\or April\or
May\or June\or
July\or August\or September\or October\or November\or December\fi}

\newcount\mins  \newcount\hours  \hours=\time \mins=\time
\def\now{\divide\hours by60 \multiply\hours by60 \advance\mins by-\hours
     \divide\hours by60         % NOTE: \divide only gives integer answers.
     \ifnum\hours>12 \advance\hours by-12
       \number\hours:\ifnum\mins<10 0\fi\number\mins\ P.M.\else
       \number\hours:\ifnum\mins<10 0\fi\number\mins\ A.M.\fi}

%**************** PAGE HEADERS *************************

\nopagenumbers
\def\runningtitle{\smallcaps Schubert calculus}
\headline={\ifnum\pageno>1\eoheadline\else\firstheadline\fi}
\def\names{\smallcaps s.\ griffeth and a.\ ram}
\def\firstheadline{}
\def\eoheadline{\ifodd\pageno\oddheadline\else\evenheadline\fi}
\def\oddheadline{\tenrm\hfil\runningtitle\hfil\folio}
\def\evenheadline{\tenrm\folio\hfil{\names}\hfil}

%**************** TITLE *************************
\vphantom{$ $}  %My kludge to get the first page to move down a bit
\vskip.75truein
\centerline{\titlefont Affine Hecke algebras and the}
\bigskip
\centerline{\titlefont Schubert calculus}
\bigskip
\centerline{\rm Stephen Griffeth}
%${}^\ast$ 
\centerline{Department of Mathematics}
\centerline{University of Wisconsin, Madison}
\centerline{Madison, WI 53706 USA}
\centerline{{\tt griffeth@math.wisc.edu}}
\bigskip
\centerline{\rm Arun Ram}
%${}^\ast$ 
\centerline{Department of Mathematics}
\centerline{University of Wisconsin, Madison}
\centerline{Madison, WI 53706 USA}
\centerline{{\tt ram@math.wisc.edu}}
%\medskip
%\centerline{Version: \today}
\bigskip\bigskip
\centerline{\sl Dedicated to Alain Lascoux}

\footnote{}{\tinyrm 
%${}^\ast$ 
\hskip-.3in
Research partially supported by the National Science Foundation
(DMS-0097977) and the National Security Agency (MDA904-01-1-0032).
}

\bigskip

%**************** ABSTRACT *************************
%\noindent{\bf Abstract.} 

%\bigskip

\section 0. Introduction

\bigskip
Using a combinatorial approach which avoids geometry,
this paper studies the ring structure of $K_T(G/B)$, the
$T$-equivariant $K$-theory of the (generalized) flag variety $G/B$.
Here, the data $G\supseteq B\supseteq T$ is a
complex reductive algebraic group (or symmetrizable 
Kac-Moody group) $G$, a Borel subgroup $B$,
and a maximal torus $T$, and $K_T(G/B)$ is the Grothendieck
group of $T$-equivariant coherent sheaves on $G/B$.
Because of the $T$-equivariance the ring $K_T(G/B)$ is an
$R$-algebra, where $R$ is the representation ring of $T$.
As explained by Grothendieck [Gd] (in the non Kac-Moody case)
and Kostant and Kumar [KK] (in the general Kac-Moody case), 
the ring $K_T(G/B)$ has 
a natural $R$-basis $\{[{\cal O}_{X_w}]\ |\ w\in W\}$, where $W$ is 
the Weyl group and ${\cal O}_{X_w}$ is the structure sheaf of the 
Schubert variety $X_w\subseteq G/B$.  One of the main problems in the
field is to understand the structure constants of the ring
$K_T(G/B)$ with this basis, that is, the coeffients $c_{wv}^z$ in the
equations
$$[{\cal O}_{X_w}][{\cal O}_{X_v}]
= \sum_{z\in W} c_{wv}^z [{\cal O}_{X_z}].
\formula$$
Our approach is to work completely combinatorially
and define $K_T(G/B)$ as a quotient of the affine nil-Hecke algebra.  
The fact that
the combinatorial approach coincides with the geometric one is a 
consequence of the results of Kostant and Kumar [KK] and Demazure [D].
In the combinatorial literature the elements $[{\cal O}_{X_w}]$
are often called (double) Grothendieck polynomials.  

Let $P$ be the weight lattice of $G$ and, for $\lambda\in P$,
let $[X^\lambda]$ be the homogeneous line bundle on $G/B$
corresponding to the character of $T$ indexed by $\lambda$.  
The theorem of Pittie [P] 
%(also proved by V.\ Snaith and R.\ Seymour)
says that the ring $K_T(G/B)$ is generated by the 
$[X^\lambda]$, $\lambda\in P$.  Steinberg [St] strengthened this
result by displaying specific $[X^{-\lambda_w}]$, $w\in W$, which
form an $R$-basis of $K_T(G/B)$.  These results are often collectively 
known as the ``Pittie-Steinberg theorem''.  

The theorems which we prove in Section 2 are simply different points of 
view on the Pittie-Steinberg theorem.  Though we are not aware of
any reference which states these theorems in the generality
which we consider, these theorems should be considered
well known.

Let $s_1,\ldots, s_n$ be the simple reflections in $W$ (determined
by the data $(G\supseteq B\supseteq T)$), let $w_0$ be the longest
element of $W$ and let $P^+$ be the set of dominant weights in $P$.
The Schubert varieties $X_{w_0s_i}$ are the codimension one Schubert
varieties in $G/B$.
In section 3 we prove ``Pieri-Chevalley'' formulas for the products
$$[X^\lambda][{\cal O}_{X_w}],
\qquad
[X^{-\lambda}][{\cal O}_{X_w}],
\qquad
[X^{w_0\lambda}][{\cal O}_{X_w}],
\qquad\hbox{and}\qquad
[{\cal O}_{X_{w_0s_i}}][{\cal O}_{X_w}],
\formula
$$
for $\lambda\in P^+$, $w\in W$ and $1\le i\le n$.
All of these Pieri-Chevalley formulas are given in terms of the combinatorics
of the Littelmann path model [Li1-3].
The formula which we give for the first product in (0.2) is due to 
Pittie and Ram [PR1].  In this paper we provide more details
of proof than appeared in [PR1].  The other formulas for the products in (0.2)
follow by applying the duality theorem of Brion [Br, Theorem 4]
to the first formula.  However, here we give an independent, 
combinatorial, proof and deduce Brion's result as a consequence.  
The last formula is a consequence of the nice formula
$$[{\cal O}_{X_{w_0s_i}}] = 1-e^{w_0\omega_i}[X^{-\omega_i}],
\formula$$
which is an easy consequence of the first two Pieri-Chevalley rules.
\smallskip\noindent

It is not difficult to ``specialize'' product formulas for 
$K_T(G/B)$ to corresponding product formulas for
$K(G/B)$, $H_T^*(G/B)$, and $H^*(G/B)$ (by using the Chern character
and comparing lowest degree terms, and ignoring the $T$-action).  
Thus the products which are computed in this paper also give results for
ordinary Grothendieck polynomials, double Schubert polynomials,
and ordinary Schubert polynomials.  In section 4 we explain how to do 
these conversions.  For most of these cases
the specialized versions of our Pieri-Chevalley rules are already very well known
(see, for example, [Ch]).

In Section 5 we give explicitly
\smallskip\noindent
\itemitem{(a)} two different kinds of formulas for $[{\cal O}_{X_w}]$
in terms of $X^\lambda$, and
\smallskip\noindent
\itemitem{(b)} complete computations of the products in (0.1)
\smallskip\noindent
for the rank two root systems.
This data allows us
to make a ``positivity conjecture'' for the coefficients
$c_{wv}^z$ in (0.1).
This conjecture generalizes the theorems of Brion [Br, formula before
Theorem 1] 
and Graham [Gr, Corollary 4.1], which treat the cases $K(G/B)$ and 
$H_T^*(G/B)$, respectively.

\medskip\noindent
{\bf Acknowledgement.}  It is a pleasure to thank Alain Lascoux
for setting the foundations of the subject of this paper.
Our approach is heavily influenced by his teachings.  In particular,
he has always promoted the study of the flag variety by divided 
difference operators (the affine, or graded, nil-Hecke algebra),
it is his work with Fulton in [FL] that provided the motivation
for the Pieri-Chevalley rules as we present them, and it his idea
of ``transition'' (see, for example, the beautiful paper [La])
which allows us to obtain product formulas
for Schubert classes in the form which we have given
in Section 5 of this paper.

\vfill\eject

\section 1. Preliminaries

\bigskip
Fix the following data and notation:
$$\matrix{
\fh^* \hfill 
&&\hbox{is a real vector space of dimension $n$}, \hfill \cr
R\hfill 
&\quad&\hbox{is a reduced irreducible root system in $\fh^*$},\hfill \cr
R^+\hfill 
&&\hbox{is a set of positive roots in $R$},\hfill \cr
W\hfill 
&&\hbox{is the Weyl group of $R$},\hfill \cr
s_1,\ldots, s_n\hfill 
&&\hbox{are the simple reflections in $W$},\hfill \cr
m_{ij}\hfill 
&&\hbox{is the order of $s_is_j$ in $W$, $i\ne j$,}\hfill \cr
R(w) = \{ \alpha\in R^+\ |\ w\alpha\not\in R^+\} \hfill
&&\hbox{is the inversion set of $w\in W$,}\hfill \cr
\ell(w)={\rm Card}(R(w)) \hfill 
&&\hbox{is the length of $w\in W$,}\hfill \cr
\le \hfill 
&&\hbox{is the Bruhat-Chevalley order on $W$,}\hfill \cr
\alpha_1,\ldots,\alpha_n\hfill 
&&\hbox{are the simple roots in $R^+$},\hfill \cr
\omega_1,\ldots, \omega_n\hfill 
&&\hbox{are the fundamental weights},\hfill \cr
P = \sum_{i=1}^n \ZZ\omega_i\hfill 
&&\hbox{is the weight lattice},\hfill \cr
P^+ = \sum_{i=1}^n \ZZ_{\ge0}\omega_i\hfill 
&&\hbox{is the set of dominant integral weights}.\hfill \cr
}$$
For a brief, easy, introduction to root systems with lots
of pictures for visualization see [NR].  
By [Bou VI \S 1 no.\ 6 Cor.\ 2 to Prop.\ 17],
if $w=s_{i_1}\cdots s_{i_p}$ be a reduced word for $w$,
then
$$
R(w) = \{\alpha_{i_p}, s_{i_p}\alpha_{i_{p-1}},\ldots, 
s_{i_p}\cdots s_{i_2}\alpha_{i_1}\},
\formula
$$ 

The {\it affine nil-Hecke algebra} is the algebra $\tilde H$ given by
generators $T_1,\ldots, T_n$ and $X^\lambda$, $\lambda\in P$,
with relations
$$
T_i^2 = T_i, \qquad
\underbrace{T_iT_jT_i\cdots}_{m_{ij}\ {\rm factors}}
=
\underbrace{T_jT_iT_j\cdots}_{m_{ij}\ {\rm factors}}, 
\qquad
X^\lambda X^\mu = X^{\lambda+\mu},
\formula
$$
and
$$
X^\lambda T_i = T_i X^{s_i\lambda} 
+ {X^\lambda-X^{s_i\lambda}\over 1-X^{-\alpha_i} }.
\formula 
$$
Let $T_w=T_{i_1}\cdots T_{i_p}$ for a reduced word $w=s_{i_1}\cdots s_{i_p}$.
Then
$$\{ X^\lambda T_w\ |\ w\in W, \lambda\in P\}
\qquad\hbox{and}\qquad
\{ T_w X^\lambda \ |\ w\in W, \lambda\in P\}
\formula
$$
are bases of $\tilde H$.

Both the {\it nil-Hecke algebra},
$$
H = \hbox{$\ZZ$-span}\{ T_w\ |\ w\in W\},
\qquad\hbox{and}\qquad
\ZZ[X] = \hbox{$\ZZ$-span}\{ X^\lambda\ |\ \lambda\in P\}
\formula$$
are subalgebras of $\tilde H$.  The action of $W$ on $\ZZ[X]$ is
given by defining
$$wX^\lambda = X^{w\lambda}, \qquad\hbox{for $w\in W$, $\lambda\in P$,}
\formula
$$
and extending linearly.  The proof of the following theorem 
is given in [R, Theorem 1.13 and Theorem 1.17].
The first statement of the theorem is due to
Bernstein, Zelevinsky, and Lusztig [Lu, 8.1] and the second statement
is due to Steinberg [St] and is known as the Pittie-Steinberg 
theorem.

\thm  Define
$$\lambda_w = w^{-1}\sum_{s_iw<w} \omega_i,
\qquad\hbox{for $w\in W$}.
\formula$$
The center of $\tilde H$ is $Z(\tilde H) = \ZZ[X]^W$ and
each element $f\in \ZZ[X]$ has a unique expansion
$$f = \sum_{w\in W} f_wX^{-\lambda_w}, 
\qquad\hbox{with $f_w \in \ZZ[X]^W$.}
\formula
$$
\endthm

Let $\varepsilon_i = 1-T_i$ and let 
$\varepsilon_w = \varepsilon_{i_1}\cdots \varepsilon_{i_p}$ for a reduced
word $w=s_{i_1}\cdots s_{i_p}$.  Then $\varepsilon_w$ is well defined
and independent of the reduced word for $w$ since
$$
\varepsilon_i^2 = \varepsilon_i,
\qquad\hbox{and}\qquad 
\underbrace{\varepsilon_i\varepsilon_j\varepsilon_i\cdots}_{
m_{ij}\ {\rm factors}}
=
\underbrace{\varepsilon_j\varepsilon_i\varepsilon_j\cdots}_{
m_{ij}\ {\rm factors}}.
\formula
$$
The second equality is a consequence of the formulas
$$\varepsilon_w
= \sum_{v\le w} (-1)^{\ell(v)}T_v
\qquad\hbox{and}\qquad
T_w = \sum_{v\le w} (-1)^{\ell(v)}\varepsilon_v
\formula
$$
which are straightforward to verify by induction on the length of $w$.

\bigskip

\section 2. The ring $K_T(G/B)$

\bigskip
Let $H$ and $\ZZ[X]$ be as in (1.5).
The {\it trivial representation}
of $H$ is defined by the homomorphism 
${\bf 1}\colon H\to \ZZ$ given by ${\bf 1}(T_i)=1$.
The first of the maps 
$$\matrix{
\ZZ[X] &\mapright{\sim} &\tilde H T_{w_0} &\mapright{\sim} 
&\tilde H\otimes_{H} {\bf 1} \cr
f &\longmapsto &fT_{w_0} &\longmapsto &f\otimes {\bf 1} \cr
}$$
is an $\tilde H$-module isomorphism
if the action of $\tilde H$ on $\ZZ[X]$ is given by
$$T_i\cdot f = {X^{\alpha_i}f - s_if\over X^{\alpha_i}-1},
\qquad\hbox{for $f\in \ZZ[X]$}.
\formula
$$

The group algebra of $P$ is 
$$R = \hbox{$\ZZ$-span}\{ e^\lambda\ |\ \lambda\in P\}
\qquad\hbox{with}\qquad
e^\lambda e^\mu = e^{\lambda+\mu},
\formula$$
for $\lambda,\mu\in P$.
Extend coefficients to $R$ so that
$\tilde H_R = R\otimes_\ZZ \tilde H$
and
$R[X] = R\otimes_\ZZ \ZZ[X]$
are $R$-algebras.
Define $K_T(G/B)$ to be the $\tilde H_R$-module 
$$K_T(G/B) = \hbox{$R$-span}\{[{\cal O}_{X_w}]\ |\ w\in W\},
\formula$$
so that the $[{\cal O}_{X_w}]$, $w\in W$, are an
$R$-basis of $K_T(G/B)$, with
$\tilde H_R$-action given by
$$X^\lambda[{\cal O}_{X_1}] = e^\lambda[{\cal O}_{X_1}],
\qquad\hbox{and}\qquad
T_i[{\cal O}_{X_w}]=\cases{
[{\cal O}_{X_{ws_i}}], &if $ws_i>w$, \cr
[{\cal O}_{X_w}], &if $ws_i<w$. \cr
}
\formula$$
If $R$ is an $R[X]$-module via the $R$-algebra homomorphism given by
$$\matrix{
e\colon &R[X] &\longrightarrow &R \cr
&X^\lambda &\longmapsto &e^\lambda \cr
}
\formula
$$
then, as $\tilde H_R$-modules,
$K_T(G/B)\cong \tilde H_R\otimes_{R[X]} R_e$, where
$R_e$ is the $R$-rank 1 $R[X]$-module determined by the homomorphism $e$.

Let $Q$ be the field of fractions of $R$ and let $\overline{Q}$ be the
algebraic closure of $Q$.  For $w\in W$ let 
$$b_w\quad\hbox{in $\overline{Q}\otimes_R K_T(G/B)$\quad be determined by}
\qquad
X^\lambda b_w = e^{w\lambda} b_w,
\quad\hbox{for $\lambda\in P$.}
\formula
$$
If the $b_w$ exist, then they are a $\overline{Q}$-basis of 
$\overline{Q}\otimes_R K_T(G/B)$ since they are eigenvectors with 
distinct eigenvalues.
If $\tau_i$, $1\le i\le n$, are the operators
on $\bar Q\otimes_R K_T(G/B)$ given by
$$\tau_i = T_i - {1\over 1-X^{-\alpha_i}},
\qquad\hbox{then}\qquad
b_1 = [{\cal O}_{X_1}]
\quad\hbox{and}\quad
\tau_i b_w = b_{ws_i},\quad
\hbox{for $ws_i>w$,}\formula$$
because, a direct computation with relation (1.3) gives that 
$
X^\lambda\tau_i b_w 
%&= X^\lambda\left( T_i - {1\over 1-X^{-\alpha_i}}\right) b_w 
%= \left(
%T_iX^{s_i\lambda}
%+{X^\lambda-X^{s_i\lambda}\over 1-X^{-\alpha_i} }
%- {X^\lambda\over 1-X^{-\alpha_i}}\right) b_w  \cr
= \tau_i X^{s_i\lambda} b_w
= \tau_i e^{ws_i\lambda} b_w 
= e^{ws_i\lambda} b_{ws_i}. 
$
Thus the $b_w$, $w\in W$, exist and the form of the 
$\tau$-operators shows that, in fact, they form
a $Q$-basis of $Q\otimes_R K_T(G/B)$ (it was
not really necessary to extend coefficients all the way to $\overline{Q}$).  
Equations (2.6) and (2.7) force
$$
\underbrace{\tau_i\tau_j\tau_i\cdots}_{m_{ij}\ {\rm factors}}
=
\underbrace{\tau_j\tau_i\tau_j\cdots}_{m_{ij}\ {\rm factors}}, \qquad
\qquad\hbox{and the equality}\qquad
\tau_i^2 = {1\over (X^{\alpha_i}-1)(X^{-\alpha_i}-1) }
$$
is checked by direct computation using (1.3).
Let $\tau_w = \tau_{i_1}\cdots \tau_{i_p}$ for a reduced word
$w=s_{i_1}\cdots s_{i_p}$.  Then, for $w\in W$,
$$b_w = \tau_{w^{-1}}b_1,
\qquad
[{\cal O}_{X_w}] = T_{w^{-1}}[{\cal O}_{X_1}] 
\qquad\hbox{and we define}\qquad
[{\cal I}_{X_w}] = \varepsilon_{w^{-1}}[{\cal O}_{X_1}], 
\formula
$$
where $\varepsilon_w$ is as in (1.11).  In terms of
geometry, $[{\cal O}_{X_w}]$ is the class of the structure sheaf
of the Schubert variety $X_w$ in $G/B$ and, up to a sign,
$[{\cal I}_{X_w}]$ is class of the sheaf ${\cal I}_{X_w}$ determined 
by the exact sequence $0\to {\cal I}_{X_w} 
\to {\cal O}_{X_w}\to {\cal O}_{\partial X_w}\to 0$, where
$\partial X_w = \bigsqcup_{v<w} BvB$ (see [Ma, Theorem 2.1(ii)] and 
[LS, equation (4)].  We are not aware of a good geometric characterization
of the basis $\{ [X^{-\lambda_w}]\ |\ w\in W\}$ of $K_T(G/B)$ which
appears in the following theorem.

\thm Let $\lambda_w$, $w\in W$, be as defined in Theorem 2.9
and let $[X^\lambda] = X^\lambda [{\cal O}_{X_{w_0}}] 
=X^\lambda T_{w_0}[{\cal O}_{X_1}]$ for 
$\lambda\in P$.  Then the $[X^{-\lambda_w}]$, 
$w\in W$, form an $R$-basis of $K_T(G/B)$.
\pf
Up to constant multiples, $[{\cal O}_{X_{w_0}}] = T_{w_0} [{\cal O}_{X_1}]$ 
is determined by the property
$$T_i[{\cal O}_{X_{w_0}}] = [{\cal O}_{X_{w_0}}],
\qquad\hbox{for all $1\le i\le n$.}
\formula$$
If constants $c_w\in Q$ are given by
$$[{\cal O}_{X_{w_0}}] = \sum_{w\in W} c_wb_w,$$
then comparing coefficients of $b_{ws_i}$, for $ws_i>w$, on each side
of (2.10) yields a recurrence relation for the $c_w$,
$$c_w=c_{ws_i} \left({1\over 1-e^{-w\alpha_i} }\right)
\quad\hbox{for $ws_i>w$,}
\qquad\hbox{which implies}\qquad
c_{w_0v^{-1}} = \prod_{\alpha\in R(v)} {1\over 1-e^{w_0\alpha}},
\formula
$$
via (1.1) and the fact that $c_{w_0}=1$.
%$$\eqalign{
%\sum_{w\in W} c_w b_w 
%&=[{\cal O}_{X_{w_0}}] 
%=T_i[{\cal O}_{X_{w_0}}] 
%= \sum_{w\in W} c_wT_i b_w \cr
%&= \sum_{ws_i>w} c_w\left(\tau_i+{1\over 1-X^{-\alpha_i}}\right) b_w
%+ c_{ws_i}\left(\tau_i+{1\over 1-X^{-\alpha_i}}\right) b_{ws_i} \cr
%&= \sum_{ws_i>w} c_wb_{ws_i}+{c_w\over 1-e^{-w\alpha_i}} b_w
%+ c_{ws_i}\tau_i^2 b_w +c_{ws_i}{1\over 1-X^{-\alpha_i}}b_{ws_i}
%\cr
%&= \sum_{ws_i>w} c_wb_{ws_i}+{c_w\over 1-e^{-w\alpha_i}} b_w
%+ {c_{ws_i}\over (e^{w\alpha_i}-1)(e^{-w\alpha_i}-1)} b_w 
%+{c_{ws_i}\over 1-e^{w\alpha_i} } b_{ws_i}.
%\cr
%}$$
%Hence $c_{ws_i} = c_w+c_{ws_i}(1-e^{w\alpha_i})^{-1}$ and
%thus
%$$c_w = c_{ws_i} \left({-e^{w\alpha_i}\over 1-e^{w\alpha_i}}\right)
%=c_{ws_i} \left({1\over 1-e^{-w\alpha_i} }\right)
%, 
%\quad\hbox{for $ws_i>w$.}
%$$
%If $w_0w=s_{i_1}\cdots s_{i_p}$ is a reduced word then,
%by Lemma ??? and the fact that $c_{w_0}=1$,
%$$\eqalign{
%c_w &= c_{w_0(w_0w)} = c_{w_0s_{i_1}\cdots s_{i_p}}
%=c_{w_0}
%\left({1\over 1-e^{-w_0s_{i_1}\cdots s_{i_{p-1}}\alpha_{i_p}} }\right)
%\cdots
%\left({1\over 1-e^{-w_0\alpha_{i_1}} }\right) \cr
%&=1\cdot \prod_{\alpha\in R(w_0w)} {1\over 1-e^{-w_0\alpha}}
%= \prod_{\alpha\in R(w)} {1\over 1-e^{-\alpha} }. \cr
%}$$
%$$
%[{\cal O}_{X_{w_0}}] = \sum_{w\in W} c_wb_w,
%\qquad\hbox{with}\qquad
%c_w = \prod_{\alpha\in R(w)} {1\over 1-e^{-\alpha}}.
%\formula$$
Thus,
$$[X^{-\lambda_v}]
= X^{-\lambda_v} [{\cal O}_{X_{w_0}}]
=\sum_{w\in W} c_w e^{-w\lambda_v} b_w,
$$
and if $C$, $M$ and $A$ are the $|W|\times |W|$ matrices given by
$$C = {\rm diag}(c_w), \quad
M = (e^{-w\lambda_v}),
\quad\hbox{and}\quad
A = (a_{zw}),
\qquad\hbox{where}\quad
b_w = \sum_{z\in W} a_{zw}[{\cal O}_{X_z}],
$$
then the transition matrix between the $X^{-\lambda_v}$ and the
$[{\cal O}_{X_z}]$ is the product $ACM$.
By (2.8) and the definition of the $\tau_i$, 
the matrix $A$ has determinant $1$.
Using the method of Steinberg [St] and subtracting row
$e^{-s_\alpha w\lambda_v}$ from row $e^{-w\lambda_v}$ in the
matrix $M$ allows one to conclude that $\det(M)$ is divisible
by 
$$\prod_{\alpha\in R^+} (1-e^{-\alpha})^{|W|/2}
\qquad\hbox{and identifying}\quad
\prod_{w\in W} e^{-w\lambda_w} =\prod_{i=1}^n \prod_{s_iw<w} e^{-\omega_i}
=(e^{-\rho})^{|W|/2}$$
as the lowest degree term determines $\det(M)$ exactly.
Thus,
$$\det(ACM)=
1\cdot \left(
\prod_{w\in W}\prod_{\alpha\in R(w)} {1\over 1-e^{-\alpha}}\right)
\left(e^{\rho}\prod_{\alpha\in R^+} \left(1-e^{-\alpha}\right) \right)^{|W|/2}
=(e^{\rho})^{|W|/2}.$$
Since this is a unit in $R$, the transition matrix between the 
$[{\cal O}_{X_w}]$ and the $X^{-\lambda_v}$ is invertible.
\endpf

\thm  The composite map
$$\matrix{
\Phi\colon 
&R[X] &\longrightarrow &\tilde H_R T_{w_0} 
&\hookrightarrow &\tilde H_R &\longrightarrow &K_T(G/B) \cr
&f &\longmapsto &fT_{w_0} 
&&h &\longmapsto &h[{\cal O}_{X_1}] \cr
}$$
is surjective with kernel
$$\ker \Phi = \langle f-e(f)\ |\ f\in R[X]^W\rangle,$$
the ideal of the ring $R[X]$ generated by the elements
$f-e(f)$ for $f\in R[X]^W$.  Hence
$$K_T(G/B) ~\cong~ {R[X]\over 
\langle f-e(f)\ |\ f\in R[X]^W\rangle }$$
has the structure of a ring.
\pf
Since $\Phi(X^\lambda) = X^\lambda T_{w_0}[{\cal O}_{X_1}] 
=X^\lambda[{\cal O}_{X_{w_0}}]$, it follows from Theorem 2.9
that $\Phi$ surjective.  Thus $K_T(G/B)\cong R[X]/\ker\Phi$.
Let $I = \langle f-e(f)\ |\ f\in R[X]^W\rangle.$  If
$f\in R[X]^W$ then, for all $\lambda\in P$,  
$$\eqalign{
\Phi(X^\lambda(f-e(f)))
&=X^\lambda(f-e(f))T_{w_0}[{\cal O}_{X_1}]
=X^\lambda T_{w_0}(f-e(f))[{\cal O}_{X_1}] \cr
&=X^\lambda T_{w_0}(e(f)-e(f))[{\cal O}_{X_1}]
=0, \cr}$$
since $f-e(f)\in Z(\tilde H_R)$.
Thus $I\subseteq \ker\Phi$.
The ring $K_T(G/B) = R[X]/\ker\Phi$ 
is a free $R$-module of rank $|W|$ and, by Theorem 1.7,
so is $R[X]/I$.  Thus $\ker\Phi = I$. 
\endpf

\bigskip\noindent

\section 3. Pieri-Chevalley formulas

\bigskip
Recall that both
$$\{ X^\lambda T_{w^{-1}}\ |\ \lambda\in P, w\in W\}
\qquad\hbox{and}\qquad
\{ T_{z^{-1}} X^\mu\ |\ \mu\in P, z\in W\}
\qquad\hbox{are bases of $\tilde H$.}
$$
If $c_{w,\lambda}^{\mu,z}\in \ZZ$ are the entries of the transition
matrix between these two bases,
$$
X^\lambda T_{w^{-1}}
= \sum_{z\in W, \mu\in P} c_{w,\lambda}^{\mu,z} T_{z^{-1}}X^\mu,
\formula$$
then applying each side of (3.1) to $[{\cal O}_{X_1}]$ gives that 
$$[X^\lambda][{\cal O}_{X_w}]
= \sum_{z\in W,\mu\in P} c_{w,\lambda}^{\mu,z} e^\mu[{\cal O}_{X_z}]\,,
\qquad\hbox{in $K_T(G/B)$.}
$$
This is the most general form of ``Pieri-Chevalley rule''.  The problem is
to determine the coefficients $c_{w,\lambda}^{\mu,z}$.

\subsection{The path model}

A {\it path} in $\fh^*$ is a piecewise linear map 
$p\colon [0,1]\to \fh^*$ such that $p(0)=0$.  
For each $1\le i\le n$ there are {\it root operators} $e_i$ and $f_i$
(see [L3] Definitions 2.1 and 2.2) which act on the paths.
If $\lambda\in P^+$ the {\it path model} for $\lambda$ is
$$
\cT^\lambda = \{ f_{i_1}f_{i_2}\cdots f_{i_l}p_\lambda \}, 
$$
the set of all paths obtained by applying the root operators to
$p_\lambda$,
where $p_\lambda$ is the straight path from $0$
to $\lambda$, that is, $p_\lambda(t)=t\lambda$, $0\le t\le 1$. 
Each path $p$ in ${\cal T}^\lambda$ is a concatenation of
segments 
$$p=p_{w_1\lambda}^{a_1}\otimes 
p_{w_2\lambda}^{a_2}\otimes \cdots \otimes p_{w_r\lambda}^{a_r}
\qquad\hbox{with}\qquad
w_1\ge w_2\ge\cdots \ge w_r
\quad\hbox{and}\quad
a_1+a_2+\cdots+a_r =1, \formula$$
where, for $v\in W$ and $a\in (0,1]$,  
$p_{v\lambda}^a$ is a piece of length $a$ from the 
straight line path $p_{v\lambda}=vp_\lambda$.  
If $W_\lambda={\rm Stab}(\lambda)$ then the $w_j$ should be
viewed as cosets in $W/W_\lambda$ and $\ge$ denotes the order on 
$W/W_\lambda$ inherited from the Bruhat-Chevalley order on $W$.
The total length of $p$ is the same as the total length of 
$p_\lambda$ which is assumed (or normalized) to be 1.  
For $p\in {\cal T}^\lambda$ let
$$\eqalign{
p(1) &= \sum_{i=1}^r a_iw_i\lambda
\quad\hbox{be the endpoint of $p$}, \cr
\iota(p)&=w_1, \quad\hbox{the initial direction of $p$,
\quad and} \cr
\phi(p)&=w_r, \quad\hbox{the final direction of $p$}. \cr
}$$

If $h\in {\cal T}^\lambda$ is such that $e_i(h)=0$ then
$h$ is the {\it head} of its {\it $i$-string}
$$S^\lambda_i(h) = \{ h, f_i h, \ldots, f_i^m h\},$$
where $m$ is the smallest positive integer such that $f_i^mh\ne 0$
and $f_i^{m+1}h =0$.  The full path model ${\cal T}^\lambda$ is
the union of its $i$-strings.  The endpoints and the
inital and final directions of the paths in the $i$-string
$S^\lambda_i(h)$ have the following properties:
$$
\matrix{
\hfill (f_i^kh)(1) &= h(1)-k\alpha_i, 
\qquad\hbox{for $0\le k\le m$}, \hfill \cr
\cr
\hbox{either} \quad \hfill
&\iota(h)=\iota(f_ih)=\cdots=\iota(f_i^mh)<s_i\iota(h) \hfill \cr
\quad\hbox{or}\hfill 
&\iota(h)<\iota(f_ih)=\cdots=\iota(f_i^mh)=s_i\iota(h),
\qquad\hbox{and} \hfill \cr
\cr
\hbox{either} \hfill 
&s_i\phi(f_i^mh)<\phi(h)=\cdots=\phi(f_i^{m-1}h)=\phi(f_i^mh) 
\hfill \cr
\quad\hbox{or} \hfill &s_i\phi(f_i^mh)=\phi(h)=\cdots
=\phi(f_i^{m-1}h) < \phi(f_i^mh). \hfill \cr
}\formula$$
The first property is [L2] Lemma 2.1a, the second is 
is [L1] Lemma 5.3, and the last is a result of applying
[L2] Lemma 2.1e to [L1] Lemma 5.3.
All of these facts are really coming from the explicit form of the 
action of the root operators on the paths in ${\cal T}^\lambda$
which is given in [L1] Proposition 4.2. 

Let $\lambda\in P^+$, $w\in W$ and $z\in W/W_\lambda$,
and let $p\in {\cal T}^\lambda$ be such that 
$\iota(p)\le wW_\lambda$ and $\phi(p)\ge z$.
Write $p$ in the form (3.2) and let $\tilde w_1,\ldots, \tilde w_r,
\tilde z$ be the maximal (in Bruhat order) coset representatives
of the cosets $w_1,\ldots, w_r,z$ such that
$$w\ge \tilde w_1\ge \tilde w_2\ge \cdots \ge \tilde w_r\ge \tilde z.
\formula$$

\thm  Recall the notation $\varepsilon_v$ from (1.11). 
Let $\lambda\in P^+$ and let $W_\lambda={\rm Stab}(\lambda)$.
Let $w\in W$.  Then, in the affine nil-Hecke algebra $\tilde H$,
$$
X^\lambda T_{w^{-1}}
= \sum_{p\in {\cal T}^\lambda\atop \iota(p)\le wW_\lambda}
T_{\phi(p)^{-1}}X^{p(1)}
\qquad\hbox{and}\qquad
X^{\lambda}\varepsilon_{w^{-1}}
= \sum_{p\in {\cal T}^{\lambda}\atop \iota(p)=w} 
\sum_{z\in W/W_\lambda\atop z\le \phi(p)} 
(-1)^{\ell(w)+\ell(z)} \varepsilon_{\tilde z^{-1}} X^{p(1)},
$$
where, if $W_\lambda \ne \{1\}$ then
$T_{\phi(p)^{-1}}=T_{\tilde w_r^{-1}}$ and $\varepsilon_{z^{-1}}
=\varepsilon_{\tilde z^{-1}}$ 
with $\tilde w_r$ and $\tilde z$ as in (3.4). 
\pf
(a) The proof is by induction on $\ell(w)$.  Let $w=s_iv$ where
$s_iv>v$.  Define 
$${\cal T}_{\le w}^\lambda 
= \{ p\in {\cal T}^\lambda\ |\ \iota(p)\le wW_\lambda\}.$$
Assume $w=s_iv>v$.  Then the facts in (3.3) imply that 
\item{(1)} ${\cal T}^\lambda_{\le w}$ is a union of the
strings $S_i(h)$ such that $h\in {\cal T}^\lambda_{\le v}$,\quad and
\item{(2)} If $h\in {\cal T}^\lambda_{\le v}$ then
either $S_i(h)\subseteq {\cal T}^\lambda_{\le v}$ or
$S_i(h)\cap {\cal T}^\lambda_{\le v} = \{h\}$.
\smallskip\noindent
Using the facts in (3.3),
a direct computation with the relation (1.3) establishes that, if
$h\in {\cal T}^\lambda_{\le v}$ then
$$\eqalign{
\sum_{p\in S_i(h)} T_{\phi(p)^{-1}}X^{\eta(1)}
&= T_{\phi(h)^{-1}}X^{h(1)}T_i, \qquad\hbox{and}\cr
\sum_{p\in S_i(h)} T_{\phi(p)^{-1}}X^{\eta(1)}
&= \cases{
T_{\phi(h)^{-1}}X^{h(1)}T_i, 
&if $S_i(h)\subseteq \cT^\lambda_{\le v}$, \cr
T_{\phi(h)^{-1}}X^{h(1)}T_i,
&if $S_i(h)\cap \cT^\lambda_{\le v} = \{h\}$. \cr
} \cr
}$$ 
Thus
$$\eqalign{
X^\lambda T_{w^{-1}}
&= X^\lambda T_{v^{-1}}T_i 
=\left(\sum_{p\in \cT^\lambda_{\le v}} 
T_{\phi(p)^{-1}}X^{p(1)}\right) T_i
\qquad\qquad\qquad\hbox{(by induction)} 
\cr 
&= \sum_{h\in \cT^\lambda_{\le v}\atop e_i(h)=0} 
\left(\sum_{S_i(h)\subseteq \cT^\lambda_{\le v}}
\sum_{p\in S_i(h)} T_{\phi(p)^{-1}}X^{p(1)} 
+
\sum_{S_i(h)\cap \cT^\lambda_{\le v}=\{h\} }
T_{\phi(h)^{-1}}X^{h(1)} \right)T_i \cr
&= \sum_{h\in \cT^\lambda_{\le w}\atop e_i(h)=0} 
\left(\sum_{S_i(h)\subseteq \cT^\lambda_{\le v} }
T_{\phi(h)^{-1}}X^{h(1)}T_i +
\sum_{S_i(h)\cap \cT^\lambda_{\le v}=\{h\} }
T_{\phi(h)^{-1}}X^{h(1)} \right)T_i \cr
&= \sum_{h\in \cT^\lambda_{\le w}\atop e_i(h)=0} 
\left(\sum_{S_i(h)\subseteq \cT^\lambda_{\le v} }
T_{\phi(h)^{-1}}X^{h(1)}T_i +
\sum_{S_i(h)\cap \cT^\lambda_{\le v}=\{h\} }
\sum_{p\in S_i(h)} T_{\phi(p)^{-1}}X^{p(1)} \right) \cr
&= \sum_{p\in \cT^\lambda_{\le w}} T_{\phi(p)^{-1}} X^{p(1)}.
}$$
(b) The proof is similar to case (a).
For $w\in W$ let
$${\cal T}^\lambda_{=w}
=\{ p\in {\cal T}^\lambda\ |\ \iota(p)=wW_\lambda\}.
$$
Assume $w=s_iv>v$.  Then the facts in (3.3) imply that 
\item{(1)} ${\cal T}^\lambda_{=w}$ is a union of the
strings $S_i(h)$ such that $h\in {\cal T}^\lambda_{=h}$, and
\item{(2)} If $h\in {\cal T}^\lambda_{=v}$ then
either $S_i(h)\subseteq {\cal T}^\lambda_{=v}$ or
$S_i(h)\cap {\cal T}^\lambda_{=v} = \{h\}$.
\smallskip\noindent
Let
$${\cal E}_{\phi(p)} = \sum_{z\in W/W_\lambda\atop z\le \phi(p)}
(-1)^{\ell(z)}\varepsilon_{\tilde z^{-1}}.
\formula$$
\smallskip\noindent
Using (3.3), a direct computation with the relation (1.3) 
establishes that, if $h\in {\cal T}^\lambda_{=v}$ with $e_ih=0$ then
$$\sum_{p \in S_i(h)} {\cal E}_{\phi(p)}X^{p(1)}T_i=0,
\qquad\hbox{and}\qquad
{\cal E}_{\phi(h)}X^{h(1)}T_i 
= -\sum_{p\in S_i(h)-\{h\} } {\cal E}_{\phi(p)}X^{p(1)}.
$$
Thus
$$\eqalign{
X^{\lambda} \varepsilon_{w^{-1}}
&= X^{\lambda} \varepsilon_{v^{-1}}\varepsilon_i
=(-1)^{\ell(v)}
\left(\sum_{p\in {\cal T}^\lambda_{=v}} 
{\cal E}_{\phi(p)}X^{p(1)}\right) T_i \cr
&=(-1)^{\ell(v)}
\left(
\sum_{S_i(h)\subseteq \cT^{\lambda}_{=v} }
\sum_{p\in S_i(h)} {\cal E}_{\phi(p)}X^{p(1)}
+
\sum_{S_i(h)\cap \cT^{\lambda}_{=v} =\{h\}}
{\cal E}_{\phi(h)}X^{h(1)} \right)T_i \cr
&=(-1)^{\ell(v)}
\left( 0 -
\sum_{S_i(h)\cap \cT^{\lambda}_{=v} =\{h\}}
\sum_{p\in S_i(h)-\{h\} } 
{\cal E}_{\phi(p)}X^{p(1)} \right) \cr
&=(-1)^{\ell(w)}
\left( \sum_{p\in {\cal T}^\lambda_{=w} } 
{\cal E}_{\phi(p)}X^{p(1)} \right).
\qquad\hbox{\qed}
}$$

\cor Let $\lambda,\mu\in P^+$ and let $w\in W$.  Then, in the 
affine nil-Hecke algebra $\tilde H$,
$$\eqalign{
X^{-\lambda}T_{w^{-1}} &= \sum_{p\in {\cal T}^{-w_0\lambda}
\atop \phi(p)=ww_0}
\sum_{z\in W/W_{-w_0\lambda}\atop zw_0\ge \iota(p)}
(-1)^{\ell(w)+\ell(z)} T_{\tilde z^{-1}}X^{p(1)}
\qquad\hbox{and}\cr
\cr
X^{w_0\mu}T_{w^{-1}} 
&= \sum_{p\in {\cal T}^\mu \atop \phi(p)=ww_0}
\sum_{z\in W/W_\mu\atop zw_0\le \phi(p)}
(-1)^{\ell(w)+\ell(z)} T_{\tilde z^{-1}}X^{p(1)}.
}
$$
\pf
The second identity is a restatement of the first with a change of 
variable $\mu = -w_0\lambda$.
The first identity is obtained by applying the algebra involution
$$
\matrix{ \tilde H &\longrightarrow &\tilde H \cr
T_w &\longmapsto &\varepsilon_w \cr
X^\lambda &\longmapsto &X^{-\lambda} \cr
}
\qquad\hbox{and the bijection} \qquad
\matrix{ {\cal T}^\lambda &\longrightarrow &{\cal T}^{-w_0\lambda} \cr
p &\longrightarrow &p^* \cr
}
$$
where $p^*$ is the same path as $p$ except translated
so that its endpoint is at the origin.
Representation theoretically, this bijection corresponds to the 
fact that $L(\lambda)^* \cong L(-w_0\lambda)$, if $L(\lambda)$
is the simple $G$-module of highest weight $\lambda$.
Note that $p^*(1) = -p(1)$, $\iota(p^*)=\phi(p)w_0$, and
$\phi(p^*)=\iota(p)w_0$.  
\endpf
 
Applying the identities from Theorem 3.5 and Corollary 3.7 
to $[{\cal O}_{X_1}]$ yields the following product formulas in 
$K_T(G/B)$.  In particular,
this gives a combinatorial proof of the ($T$-equivariant extension)
of the duality theorem of Brion [Br, Theorem 4].
For $\lambda\in P$ and $w\in W$ let 
$[X^\lambda]=X^\lambda[{\cal O}_{X_{w_0}}]
=X^\lambda T_{w_0}[{\cal O}_{X_1}]$
and let $c_{\lambda,w}^z$ be given by
$$[X^\lambda][{\cal O}_{X_w}]
= \sum_{z\in W} c_{\lambda,w}^z [{\cal O}_{X_z}],
\formula$$

\cor  Let $\lambda\in P^+$, $w\in W$ and $W_\lambda = {\rm Stab}(\lambda)$.
Then, with notation as in (3.8),
$$\displaylines{
c_{\lambda,w}^z
=\sum_{p\in {\cal T}^\lambda\atop 
wW_\lambda\ge \iota(p)\ge \phi(p)=zW_\lambda} e^{p(1)},
\cr
\cr
c_{w_0\lambda,w}^z
= (-1)^{\ell(w)+\ell(z)} c_{\lambda,zw_0}^{ww_0},
\qquad\hbox{and} \qquad
c_{-\lambda,w}^z
= (-1)^{\ell(w)+\ell(z)}c_{-w_0\lambda,zw_0}^{ww_0}.
}
$$
\endthm

\prop For $1\le i\le n$,
\quad $[{\cal O}_{X_{w_0s_i}}] = 1-e^{w_0\omega_i}[X^{-\omega_i}].$
\pf We shall show that 
$$X^{-\omega_i}[{\cal O}_{X_{w_0}}]
= e^{-w_0\omega_i}([{\cal O}_{X_{w_0}}]-[{\cal O}_{X_{w_0s_i}}]),
\formula$$
and the result will follow by solving for $[{\cal O}_{X_{s_iw_0}}]$.
Let $\omega_j=-w_0\omega_i$.  By Corollary 3.9,
$$c_{-\omega_i,w_0}^z
=(-1)^{\ell(w_0)+\ell(z)} c_{\omega_j,zw_0}^{1}
=(-1)^{\ell(w_0)+\ell(z)} \sum_{p\in {\cal T}^{\omega_j}
\atop zw_0\ge \iota(p)\ge \phi(p)=1} e^{p(1)}.$$
The straight line path to $\omega_j$, $p_{\omega_j}$, has
$\iota_{zw_0}(p_{\omega_j})=\phi_{zw_0}(\omega_j)$ and
is the unique path in ${\cal T}^{\omega_j}$ which may have final
direction $1$. Suppose $\phi_{zw_0}(p_{\omega_j})=1$.  Then, since
$s_j$ is the only simple reflection which is not in 
${\rm Stab}(\omega_j)$, 
it must be that $zw_0\not\ge s_k$ for all $k\ne j$.  Thus
$zw_0=1$ or $zw_0 = s_j$ and so $c_{-\omega_i,w_0}^z\ne 0$ only
if $z=w_0$ or $z = s_jw_0 = w_0s_i$.  Now (3.11) follows
since $p_{\omega_j}$ has endpoint $\omega_j = -w_0\omega_i$. 
\endpf

\cor  Let $c_{wv}^z$ be as in (3.8).  Then, for $1\le i\le n$,
$c_{w_0s_i,w}^w = -(e^{-(w\omega_i-w_0\omega_i)}-1)$, and
$$c_{w_0s_i,w}^z
=(-1)^{\ell(w)+\ell(z)+1}
\sum_{p\in {\cal T}^{-w_0\omega_i}\atop 
zw_0\ge \iota(p)\ge \phi(p)=ww_0}
e^{w_0\omega_i+p(1)},
\qquad\hbox{for $z\ne w$}.
$$
\pf
This follows from Proposition 3.10 and Corollary 3.9 and the fact
that, in the case when $z=w$, there is a unique path $p$ with
$ww_0=\iota(p)=\phi(p)=ww_0$ and endpoint
$p(1)=ww_0(-w_0\omega_i)=-w\omega_i$.
\endpf

\bigskip

\section 4. Converting to $H^*_T(G/B)$

\bigskip

The {\it graded nil-Hecke algebra} is the algebra $H_{\rm gr}$ given by
generators $t_1,\ldots, t_n$ and $x_\lambda$, $\lambda\in P$,
with relations
$$
t_i^2 = 0, \qquad
\underbrace{t_it_jt_i\cdots}_{m_{ij}\ {\rm factors}}
=
\underbrace{t_jt_it_j\cdots}_{m_{ij}\ {\rm factors}}, \qquad
x_{\lambda+\mu} = x_{\lambda}+x_\mu, 
\quad\hbox{and}\quad
x_\lambda t_i = t_i x_{s_i\lambda} 
+ \langle \lambda,\alpha_i^\vee\rangle. 
\formula
$$
The subalgebra of $H_{\rm gr}$ generated by the $x_\lambda$ is 
the polynomial ring $\ZZ[x_1,\ldots, x_n]$, where $x_i = x_{\omega_i}$,
and $W$ acts on $\ZZ[x_1,\ldots, x_n]$ by
$$wx_\lambda = x_{w\lambda}
\qquad\hbox{and}\qquad
w(fg) = (wf)(wg),
\qquad\hbox{for $w\in W$, $\lambda\in P$, 
$f,g\in \ZZ[x_1,\ldots, x_n]$.}$$
Then the last formula in (4.1) generalizes to
$$f t_i = t_i (s_if)
+ {f-s_if\over \alpha_i},
\qquad\hbox{for $f\in \ZZ[x_1,\ldots, x_n]$}.
$$
Let $t_w = t_{i_1}\cdots t_{i_p}$ for a reduced word
$w=s_{i_1}\cdots s_{i_p}$ and let $\ZZ W^*$ be the subalgebra
of $H_{\rm gr}$ spanned by the $t_w$, $w\in W$.  Then
$$
\{x_1^{m_1}\cdots x_n^{m_n} t_w\ |\ w\in W, \ \ 
m_i\in \ZZ_{\ge 0} \}
\qquad\hbox{and}\qquad
\{t_w x_1^{m_1}\cdots x_n^{m_n} \ |\ w\in W, \ \ 
m_i\in \ZZ_{\ge 0} \} 
$$
are bases of $H_{\rm gr}$.   

Let $S = \ZZ[y_1,\ldots, y_n]$ and extend coefficients to $S$
so that $H_{{\rm gr},S} = S\otimes_\ZZ H_{\rm gr}$ and 
$S[x_1,\ldots,x_n]=S\otimes_\ZZ \ZZ[x_1,\ldots, x_n]$ are $S$-algebras.
Define $H^*_T(G/B)$ to be the $H_{{\rm gr},S}$ module
$$H^*_T(G/B) = \hbox{$S$-span}\{ [X_w]\ |\ w\in W\},
\formula$$
so that the $[X_w]$, $w\in W$, are an $S$-basis of $K_T(G/B)$,
with $H_{{\rm gr}, S}$-action given by
$$x_i[X_1] = y_i[X_1],
\qquad\hbox{and}\qquad
t_i[X_w] = \cases{ 
[X_{ws_i}], &if $ws_i>w$, \cr
0, &if $ws_i<w$, \cr}
\formula
$$

Let $y$ be the $S$-algebra homomorphism given by
$$\matrix{y\colon &S[x_1,\ldots, x_n] 
&\longrightarrow &S \cr
&x_i &\longmapsto &y_i \cr
}$$
so that $H^*_T(G/B) \cong H_{{\rm gr},S}
\otimes_{S[x_1,\ldots, x_n]} y$ as $H_{{\rm gr},S}$-modules
Then, using analogous methods to the $K_T(G/B)$ case
proves the following
theorem, which gives the ring structure of $H^*T(G/B)$
(see also the proof of [KR, Prop.\ 2.9] for the same argument
with (non-nil) graded Hecke algebras). 

\thm  The composite map
$$\matrix{
\Phi\colon 
&S[x_1,\ldots, x_n] &\longrightarrow &H_{{\rm gr}, S} t_{w_0} 
&\hookrightarrow &H_{{\rm gr},S} &\longrightarrow &H^*_T(G/B) \cr
&f &\longmapsto &ft_{w_0} 
&&h &\longmapsto &h[X_1] \cr
}$$
is surjective with kernel
$$\ker \Phi = \langle f-y(f)\ |\ f\in S[x_1,\ldots, x_n]^W\rangle,$$
the ideal of the ring $S_[x_1,\ldots, x_n]$ generated by the elements
$f-y(f)$ for $f\in S[x_1,\ldots,x_n]^W$.  Hence
$$H^*_T(G/B) ~\cong~ {\ZZ[y_1,\ldots, y_n, x_1,\ldots, x_n]\over 
\langle f-y(f)\ |\ f\in S[x_1,\ldots, x_n]^W\rangle }$$
has the structure of a ring.
\endthm

As a vector space $H_{\rm gr} = 
\ZZ[x_1,\ldots, x_n]\otimes \ZZ W_{\rm gr}$.
Let $\widehat{ H_{\rm gr} }= \QQ[[x_1,\ldots, x_n]]\otimes \QQ W_{gr}$
with multiplication determined by the relations in (4.1).  Then
$\widehat{ H_{\rm gr} }$ is a completion of $H_{\rm gr}$ (this simply
allows us to write infinite sums) 
and the elements of $\widehat{ H_{\rm gr} }$ given by
$$
{\rm ch}(X^\lambda)=
\sum_{r\ge 0 } \hbox{$1\over r!$}\,x_\lambda^r 
\qquad\hbox{and}\qquad
{\rm ch}(T_i) = t_i\cdot{x_{\alpha_i}\over 1-{\rm ch}(X^{\alpha_i})}
\formula
$$
satisfy the relations of $\tilde H$ and thus ${\rm ch}$ extends to
a ring homomorphism
${\rm ch}\colon \tilde H \longrightarrow \widehat{H_{\rm gr} }.$
It is this fact that really makes possible the transfer from
$K$-theory to cohomomology possible.  Though is it not difficult
to check that the elements in (3.5) satisfy the defining relations
of $\tilde H$ it is helpful to realize that these
formulas come from geometry.  As explained in [PR2], the action of
$T_i$ on $K_T(G/B)$ and the action of $t_i$ on $H_T^*(G/B)$ are,
respectively, the push-pull operators 
$\pi_i^*(\pi_i)_!$ and $\pi_i^*(\pi_i)_*$, where 
if $P_i$ is a minimal parabolic subgroup of $G$ then 
$\pi_i\colon G/P_i\to G/B$ is the natural surjection.
Then the first formula in (3.5) is the definition of the Chern character, 
and the second formula is the Grothedieck-Riemann-Roch theorem applied to 
the map $\pi_i$.  The factor $\alpha_i/(1-{\rm ch}(X^{\alpha_i}))$ is
the Todd class of the bundle of tangents along the fibers of $\pi_i$
(see [Hz, page 91]).

Then $\widehat{ H^*_T}(G/B)_{\QQ} 
= \QQ[[y_1,\ldots, y_n]]\otimes_{\ZZ[y_1,\ldots, y_n]} H^*_T(G/B)$
is the appropriate completion of $H^*_T(G/B)$
to use to transfer the ring homomorphism 
${\rm ch}\colon \tilde H_R\to \widehat{ H_{\rm gr} }$ to a ring homomorphism
$${\rm ch}\colon K_T(G/B) \longrightarrow \widehat{ H^*_T}(G/B)_{\QQ}
\qquad\hbox{by setting}\quad
{\rm ch}(h[{\cal O}_{X_1}]) = {\rm ch}(h)[X_1],
\quad\hbox{for $h\in \tilde H_R$.}
\formula$$
The ring $\widehat{ H^*_T}(G/B)_{\QQ}$ is a graded ring with
$${\rm deg}(y_i)=1
\qquad\hbox{and}\qquad {\rm deg}([X_w]) = \ell(w_0)-\ell(w),
\formula$$
$$\hbox{and,\quad for $w\in W$,}\qquad\qquad 
{\rm ch}([{\cal O}_{X_w}]) = [X_w]+ ~\hbox{higher degree terms}.
\qquad\qquad\qquad\qquad \hfill
\formula$$ 
In summary, if $e_i=e^{\omega_i}$, $X_i=X^{\omega_i}$, $y_i=y_{\omega_i}$,
$x_i = x_{\omega_i}$,
$$\eqalign{
R[X] &= \ZZ[e_1^{\pm 1},\ldots, e_n^{\pm1}, X_1^{\pm1},\ldots, X_n^{\pm1}], \cr
\ZZ[X] &= \ZZ[X_1^{\pm1},\ldots, X_n^{\pm1}], \cr
}
\quad\hbox{and}\quad
\widehat{S}[x_1,\ldots, x_n] =\QQ[[y_1,\ldots, y_n]][x_1,\ldots, x_n],$$
then there is a commutative diagram of ring homomorphisms
$$\matrix{
\displaystyle{
K_T(G/B) 
={R[X] \over \langle f-e(f)\ |\ f\in R[X]^W\rangle} }\hfill
&\mapright{{\rm ch}} 
&\displaystyle{
H^*_T(G/B)_{\QQ} 
={\widehat{S}[x_1,\ldots, x_n] 
\over \langle f-y(f)\ |\ 
f\in \widehat{S}[x_1,\ldots, x_n]^W\rangle} } \hfill \cr
\cr
\cr
\mapdown{e_i=1} &&\mapdown{y_i=0} \cr
\cr
\displaystyle{
K(G/B) 
={\ZZ[X] \over \langle f-f(1)\ |\ f\in \ZZ[X]^W\rangle}  } \hfill
&\mapright{{\rm ch}}
&\displaystyle{
H^*(G/B)_{\QQ} 
= {\QQ[x_1,\ldots, x_n] \over 
\langle f-f(0)\ |\ f\in \QQ[x_1,\ldots, x_n]^W\rangle}  }.\cr
}$$

\bigskip

\section 5.  Rank two and a positivity conjecture

\bigskip
In this section we will give explicit formulas for the
rank two root systems.  
The data supports the following positivity conjecture which
generalizes the theorems of Brion [Br, formula before Theorem 1] 
and Graham [Gr, Corollary 4.1].

\conj For $\beta\in R^+$ let $y_\beta = e^{-\beta}$ and
$a_\beta = e^{-\beta}-1$ and let $d(w) = \ell(w_0)-\ell(w)$ for
$w\in W$.  Let $c_{wv}^z$ be the structure 
constants of $K_T(G/B)$ with respect to the basis
$\{ [{\cal O}_{X_w}]\ |\ w\in W\}$ as defined in (0.1).
Then
$$c_{wv}^z = (-1)^{d(w)+d(v)-d(z)}f(\alpha,y),
\qquad\hbox{where}\quad
f(\alpha,y)\in \ZZ_{\ge 0}[\alpha_\beta,y_\beta\ |\ \beta\in R^+],$$
that is, $f(\alpha,y)$ is a polynomial in the variables
$\alpha_\beta$ and $y_\beta$, $\beta\in R^+$, which has nonnegative
integral coefficients.
\endthm

In the following, for brevity, use the following notations:
$$
\matrix{
\hbox{in $K_T(G/B)$,} \hfill
&[w] = [{\cal O}_{X_w}],  \hfill
&\alpha_{rs} = e^{-(r\alpha_1+s\alpha_2)}-1, \hfill
&\hbox{and}\quad
&y_{rs} = e^{-(r\alpha_1+s\alpha_2)}, \hfill
\cr
\hbox{in $K(G/B)$,}  \hfill
&[w] = [{\cal O}_{X_w}],  \hfill
&\alpha_{rs} = 0,  \hfill
&\hbox{and}\quad
&y_{rs} = 1, \hfill
\cr
\hbox{in $H_T^*(G/B)$,} \hfill
&[w] = [X_w],  \hfill
&\alpha_{rs} = r\alpha_1+s\alpha_2, \hfill
&\hbox{and}\quad
&y_{rs} = 1, \hfill
\cr
\hbox{in $H^*(G/B)$,} \hfill
&[w] = [X_w],  \hfill
&\alpha_{rs} = 0,  \hfill
&\hbox{and}\quad
&y_{rs} = 1,\hfill 
\cr
}$$
and in $H^*_T(G/B)$ and in $H^*(G/B)$ the terms in $\{\,\}$ brackets 
do not appear.

\bigskip\noindent
{\bf Type $A_2$.}
For the root system $R$ of type $A_2$ 
$$\matrix{
\alpha_1 = -\omega_1+2\omega_2, 
\phantom{1\over3_j}
&\lambda_1 = \rho, \hfill
&\lambda_{s_1} = \omega_2={1\over3}\alpha_1+{2\over 3}\alpha_2,\ \ \hfill 
&\lambda_{s_2s_1} = s_2\omega_2
= \phantom{-}{1\over 3}\alpha_1-{1\over 3}\alpha_2, 
\hfill \cr
\alpha_2 = \phantom{-}2\omega_1-\omega_2, 
\phantom{1\over3_j}
&\lambda_{w_0}=0, 
&\lambda_{s_2} = \omega_1={2\over 3}\alpha_1+{1\over 3}\alpha_2,\ \   \hfill
&\lambda_{s_1s_2} = s_1\omega_1 = -{1\over 3}\alpha_1+{1\over3}\alpha_2. 
\hfill\cr
}$$
Formulas for the Schubert classes in terms of homogeneous line bundles
can be given by
$$\matrix{
[s_1s_2s_1]=1,\quad \hfill &\qquad
&[1] =(1-e^{s_1\omega_1}X^{-\omega_1})[s_1] 
= (1-e^{s_2\omega_2}X^{-\omega_2})[s_2],  \hfill \cr
[s_2s_1]=1-e^{-\omega_1}X^{-\omega_2},\hfill
&&[s_1s_2] = 1-e^{-\omega_2}X^{-\omega_1} \hfill \cr
[s_1]=(1-e^{s_2\omega_2}X^{-\omega_2})[s_2s_1],\hfill
&&[s_2] = (1-e^{s_1\omega_1}X^{-\omega_1})[s_1s_2],\hfill \cr
}
$$
and
$$\eqalign{
[s_1s_2s_1] &= 1, \qquad
[s_1s_2] = 1-e^{-\omega_2}X^{-\omega_1}, \qquad
[s_2s_1] = 1-e^{-\omega_1}X^{-\omega_2}, \cr
[s_1] &=
1
-e^{-\omega_2}X^{-s_1\omega_1}
-e^{-\omega_2}X^{-\omega_1}
+e^{-2\omega_2}X^{-\omega_2},
\cr
[s_2] &=
1
-e^{-\omega_1}X^{-s_2\omega_2}
-e^{-\omega_1}X^{-\omega_2}
+e^{-2\omega_1}X^{-\omega_1},
\cr
[1] &= 
1
-e^{-\omega_2}X^{-s_1\omega_1}
-e^{-\omega_1}X^{-s_2\omega_2}
+e^{-2\omega_1}X^{-\omega_1}
+e^{-2\omega_2}X^{-\omega_2}
-e^{-\rho}X^{-\rho}.
\cr
}$$
The multiplication of the Schubert classes is given by
%Multiplication by $[1]$:
$$\eqalign{
[1]^2 
&=
-\alpha_{10}\alpha_{01}\alpha_{11}[1],
\cr
[1][s_1]
&=
\alpha_{01}\alpha_{11}[1],
\cr
[1][s_2]
&=
\alpha_{10}\alpha_{11}[1],
\cr
[1][s_1s_2]
&=
-\alpha_{11}[1],
\cr
[1][s_2s_1]
&=
-\alpha_{11}[1],
\cr
}
%}$$
%
%\bigskip
%%Multiplication by $[s_1]$:
%
%$$\eqalign{
\qquad\eqalign{
[s_1]^2
&=
\alpha_{01}\alpha_{11}[s_1],
\cr
[s_1][s_2]
&=
-\alpha_{11}[1],
\cr
[s_1][s_1s_2]
&=
y_{01}[1]
-\alpha_{01}[s_1],
\cr
[s_1][s_2s_1]
&=
-\alpha_{11}[s_1],
\cr
\cr
}
%}$$
%
%\bigskip
%%Multiplication by $[s_2]$:
%
%$$\eqalign{
\qquad\eqalign{
[s_2]^2
&=
\alpha_{01}\alpha_{11}[s_2],
\cr
[s_2][s_1s_2]
&=
-\alpha_{11}[s_2], 
\cr
[s_2][s_2s_1]
&=
y_{10}[1]
-\alpha_{10}[s_2],
\cr
\cr
\cr
}$$
%\bigskip
%Multiplication by $[s_1s_2]$:
$$\eqalign{
[s_1s_2]^2
&=
y_{01}[s_2]
-\alpha_{01}[s_1s_2],
\cr
[s_1s_2][s_2s_1]
&=\,\{\,-[1]\,\}\,+[s_1]+[s_2],
\cr
}
%}$$
%
%\bigskip
%%Multiplication by $[s_2s_1]$:
%
%$$\eqalign{
\qquad\qquad\eqalign{
[s_2s_1]^2
&= y_{10}[s_1] -\alpha_{10}[s_2s_1].
\cr
\cr
}$$

\bigskip\noindent
{\bf Type $B_2$.}
For the root system $R$ of type $B_2$ 
$$\matrix{
\alpha_1 = \phantom{-}2\omega_1-\omega_2, \hfill
&\quad 
&\lambda_1 = \rho = 2\alpha_1+{3\over2}\alpha_2, \hfill
&\quad
&\lambda_{s_1} = \omega_2 = \alpha_1+\alpha_2, \hfill 
&\quad
\cr
\alpha_2 = -2\omega_1+2\omega_2, \hfill
&&\lambda_{w_0}=0, \hfill
&&\lambda_{s_2} = \omega_1 = \alpha_1+{1\over2}\alpha_2,  \hfill
\cr
}$$
$$\matrix{
\lambda_{s_2s_1} = s_2\omega_2 =\alpha_1, \hfill
&\quad
&\lambda_{s_1s_2s_1} = s_1s_2\omega_2 =-\alpha_1, \hfill 
\cr
\lambda_{s_1s_2} = s_1\omega_1 = {1\over2}\alpha_2, \hfill
&&\lambda_{s_2s_1s_2} = s_2s_1\omega_1 = -{1\over2}\alpha_2. \hfill
}$$
Formulas for the Schubert classes in terms of homogeneous line bundles
can be given by
$$\matrix{
[s_1s_2s_1s_2]=1, \hfill
&&[1] 
=(1-e^{s_1\omega_1}X^{-\omega_1})[s_1] 
= (1-e^{s_2\omega_2}X^{-\omega_2})[s_2], \hfill \cr
[s_1s_2s_1]=1-e^{-\omega_2}X^{-\omega_2},\hfill
&&[s_2s_1s_2] = 1-e^{-\omega_1}X^{-\omega_1}, \hfill\cr
[s_2s_1]=(1-e^{-\omega_1}X^{-s_1\omega_1})[s_2s_1s_2],\hfill
&&[s_1s_2] = (1-e^{s_2s_1\omega_1}X^{-\omega_1})[s_2s_1s_2],\hfill \cr
[s_1]=(1-e^{s_2\omega_2}X^{-\omega_2})[s_2s_1],\hfill
&&[s_2] = (1-e^{s_1\omega_1}X^{-\omega_1})[s_1s_2],\hfill \cr
}$$
and
$$\eqalign{
[s_1s_2s_1s_2]&=1,
\qquad 
[s_1s_2s_1] = 1-e^{-\omega_2}X^{-\omega_2},
\qquad
[s_2s_1s_2] = 1-e^{-\omega_1}X^{-\omega_1},
\cr
[s_1s_2] &=
(1-e^{-\omega_2})
-e^{-\omega_2}X^{-\omega_2}
-e^{-\omega_2}X^{-s_2\omega_2}
+(e^{-\rho}+e^{-s_1\rho})X^{-\omega_1},
\cr
[s_2s_1] &=
1
-e^{-\omega_1}X^{-\omega_1}
-e^{-\omega_1}X^{-s_1\omega_1}
+e^{-2\omega_1}X^{-\omega_2},
\cr
[s_1] &=
(1-e^{-\omega_2})
+(e^{-\rho}+e^{-s_1\rho})X^{-s_1\omega_1}
+(e^{-\rho}+e^{-s_1\rho})X^{-\omega_1} \cr
&\qquad
-e^{-\omega_2}X^{-s_1s_2\omega_2}
-e^{-\omega_2}X^{-s_2\omega_2}
-(e^{-2\omega_2}+e^{-\omega_2})X^{-\omega_2},
\cr
[s_2] &=
(1+e^{-2\omega_1})
+e^{-2\omega_1}X^{-s_2\omega_2}
+e^{-2\omega_1}X^{-\omega_2} \cr
&\qquad
-e^{-\omega_1}X^{-s_2s_1\omega_1}
-e^{-\omega_1}X^{-s_1\omega_1}
-(e^{-3\omega_1}+e^{-\omega_1})X^{-\omega_1},
\cr
[1] &=
(1+e^{-2\omega_1})
-e^{-\omega_1}X^{-s_2s_1\omega_1}
+(e^{-\rho}+e^{-s_1\rho})X^{-s_1\omega_1}
-(e^{-3\omega_1}+e^{-\omega_1})X^{-\omega_1} \cr
&\qquad\qquad\qquad
-e^{-\omega_2}X^{-s_1s_2\omega_2}
+e^{-2\omega_1}X^{-s_2\omega_2}
-(e^{-2\omega_2}+e^{-\omega_2})X^{-\omega_2}
+e^{-\rho}X^{-\rho}.
}$$
The multiplication of the Schubert classes is given by
%Multiplication by [1]
$$\eqalign{
[1]^2 
&=
\alpha_{10}\alpha_{01}\alpha_{11}\alpha_{21}[1],
\cr
[1][s_1]
&=
-\alpha_{01}\alpha_{11}\alpha_{21}[1],
\cr
[1][s_2]
&=
-\alpha_{10}\alpha_{11}\alpha_{21}[1],
\cr
[1][s_1s_2]
&=
\alpha_{11}\alpha_{21}[1],
\cr
[1][s_2s_1]
&=
\alpha_{11}\alpha_{21}[1],
\cr
[1][s_1s_2s_1]
&=
-\alpha_{11}(1+y_{11})[1], \cr
[1][s_2s_1s_2]
&=
-\alpha_{21}[1], 
}
%}$$
%
%\bigskip
%%Multiplication by $[s_1s_2s_1]$:
%
%$$\eqalign{
\qquad
\eqalign{
[s_1s_2s_1]^2
&=
\{\,-y_{11}[s_1]\,\}
+(y_{01}+y_{11})[s_2s_1]
-\alpha_{01}[s_1s_2s_1],
\cr
[s_1s_2s_1][s_2&s_1s_2]
=
\{\,[1]-[s_1]-[s_2]\,\}+[s_1s_2]+[s_2s_1],
\cr
\cr
%}$$
%
%\bigskip
%%Multiplication by $[s_2s_1s_2]$:
%
%$$\eqalign{
[s_2s_1s_2]^2
&=
y_{10}[s_1s_2] -\alpha_{10}[s_2s_1s_2],
\cr
[s_2s_1]^2
&=
-\alpha_{21}y_{10}[s_1]
+\alpha_{10}\alpha_{21}[s_2s_1],
\cr
\cr
[s_2s_1][s_1s_2s_1]
&=
y_{21}[s_1]
-\alpha_{21}[s_2s_1],
\cr
[s_2s_1][s_2s_1s_2]
&=
\{\,-y_{10}[1]\,\}
+y_{10}[s_1]
+y_{10}[s_2]
-\alpha_{10}[s_2s_1],
\cr
}$$
%Multiplication by $[s_1]$:
$$\eqalign{
[s_1]^2
&=
-\alpha_{01}\alpha_{11}\alpha_{21}[s_1],
\cr
[s_1][s_2]
&=
\alpha_{11}\alpha_{21}[1],
\cr
[s_1][s_1s_2]
&=
-\alpha_{11}(y_{01}+y_{11})[1]
+\alpha_{01}\alpha_{11}[s_1],
\cr
[s_1][s_2s_1]
&=
\alpha_{11}\alpha_{21}[s_1],
\cr
[s_1][s_1s_2s_1]
&=
-\alpha_{11}(1+y_{11})[s_1], 
\cr
[s_1][s_2s_1s_2]
&=
y_{11}[1]
-\alpha_{11}[s_1], 
\cr
}
%}$$
%
%
%\bigskip
%%Multiplication by $[s_2]$:
%
%$$\eqalign{
\qquad
\eqalign{
[s_2]^2
&=
-\alpha_{10}\alpha_{11}\alpha_{21}[s_2],
\cr
[s_2][s_1s_2]
&=
\alpha_{11}\alpha_{21}[s_2], 
\cr
[s_2][s_2s_1]
&=
-\alpha_{21}y_{10}[1]
+\alpha_{10}\alpha_{21}[s_2],
\cr
[s_2][s_1s_2s_1]
&=
y_{21}[1]
-\alpha_{21}[s_2], 
\cr
[s_2][s_2s_1s_2]
&=
-\alpha_{21}[s_2],
\cr
}$$
%
%\bigskip
%Multiplication by $[s_1s_2]$:
%
$$\eqalign{
[s_1s_2]^2
&=
-\alpha_{11}(y_{01}+y_{11})[s_2]
+\alpha_{01}\alpha_{11}[s_1s_2],
\cr
[s_1s_2][s_2s_1]
&=
(\{\,\alpha_{11}\,\}+y_{21})[1]
-\alpha_{11}[s_1]
-\alpha_{21}[s_2],
\cr
[s_1s_2][s_1s_2s_1]
&=
\{\,-(y_{01}+y_{11}) [1]\,\}
+y_{01}[s_1]
+(y_{11}+y_{12})[s_2]
-\alpha_{01}[s_1s_2], 
\cr
[s_1s_2][s_2s_1s_2]
&=
y_{11}[s_2]
-\alpha_{11}[s_1s_2],
\cr
}$$
%
%\bigskip
%Multiplication by $[s_2s_1]$:
%
$$\eqalign{
[s_2s_1]^2
&=
-\alpha_{21}y_{10}[s_1]
+\alpha_{10}\alpha_{21}[s_2s_1],
\cr
[s_2s_1][s_1s_2s_1]
&=
y_{21}[s_1]
-\alpha_{21}[s_2s_1],
\cr
[s_2s_1][s_2s_1s_2]
&=
\{\,-y_{10}[1]\,\}
+y_{10}[s_1]
+y_{10}[s_2]
-\alpha_{10}[s_2s_1],
\cr
}$$

\bigskip\noindent
{\bf Type $G_2$.}
For the root system $R$ of type $G_2$ 
$$\matrix{
\lambda_1 = \rho = 5\alpha+3\alpha_2, \hfill
&\quad
&\lambda_{s_1s_2s_1} = s_1s_2\omega_2 =\alpha_2, \hfill  
\cr
\lambda_{s_1} = \omega_2 = 3\alpha_1+2\alpha_2,\hfill 
&&\lambda_{s_2s_1s_2s_1} = s_2s_1s_2\omega_2 = -\alpha_2, \hfill 
\cr
\lambda_{s_2} = \omega_1 =2\alpha_1+\alpha_2, \hfill
&&\lambda_{s_1s_2s_1s_2} = s_1s_2s_1\omega_1 =-\alpha_1, \hfill
\cr
\lambda_{s_2s_1} = s_2\omega_2 = 3\alpha_1+\alpha_2,\hfill
&&\lambda_{s_1s_2s_1s_2s_1} = s_1s_2s_1s_2\omega_2 
= -3\alpha_1-\alpha_2, \hfill 
\cr
\lambda_{s_1s_2} = s_1\omega_1 = \alpha_1 + \alpha_2, \hfill
&&\lambda_{s_2s_1s_2s_1s_2} = s_2s_1s_2s_1\omega_1 = -\alpha_1-\alpha_2,\hfill
\cr
\lambda_{s_2s_1s_2} = s_2s_1\omega_1 = \alpha_1,\hfill
&&\lambda_{w_0}=0.\hfill
\cr
}$$
Formulas for the Schubert classes in terms of homogeneous line bundles
can be given by
$$\matrix{
[s_1s_2s_1s_2s_1s_2]=1, \hfill
&&
[1] 
=(1-e^{s_1\omega_1}X^{-\omega_1})[s_1] 
= (1-e^{s_2\omega_2}X^{-\omega_2})[s_2], \hfill \cr
[s_1s_2s_1s_2s_1]=1-e^{-\omega_2}X^{-\omega_2},\hfill
&&[s_2s_1s_2s_1s_2] = 1-e^{-\omega_1}X^{-\omega_1}, \hfill\cr
[s_2s_1s_2s_1]
=(1-e^{-\omega_1}X^{-s_1\omega_1})[s_2s_1s_2s_1s_2],\hfill
&&[s_1s_2s_1s_2] = 
(1-e^{-s_1\omega_1}X^{-\omega_1})[s_2s_1s_2s_1s_2],\hfill \cr
[s_1s_2s_1]=\hbox{see below}, \hfill
&&\displaystyle{
[s_2s_1s_2] 
= {1-e^{-s_2s_2\omega_1}X^{-\omega_1}\over 
1+X^{-\omega_1} }\,[s_1s_2s_1s_2] },\hfill \cr
[s_2s_1]=(1-e^{-\omega_1}X^{-s_1s_2s_1\omega_1})[s_2s_1s_2], \hfill
&&
[s_1s_2] 
= (1-e^{s_2s_1\omega_1}X^{-\omega_1})[s_1s_2],\hfill \cr
[s_1]=(1-e^{s_2\omega_2}X^{-\omega_2})[s_2s_1],\hfill
&&
[s_2] 
= (1-e^{s_1\omega_1}X^{-\omega_1})[s_1s_2],\hfill \cr
}$$
\smallskip
$$
[s_1s_2s_1] =
{(1-e^{-\alpha_2}X^{-\omega_2})[s_2s_1s_2s_1]
+e^{-\alpha_2}(1+e^{\omega_1}X^{-\omega_2})[s_2s_1]
\over 1+e^{-\alpha_2}
},
$$
and
$$\eqalign{
[w_0] &= 1,
\qquad
[s_2s_1s_2s_1s_2] =
1-y_{21}X^{-\omega_1},
\qquad 
[s_1s_2s_1s_2s_1] =
1-y_{32}X^{-\omega_2},
\cr
[s_2s_1s_2s_1] &=
1-y_{21}X^{-\omega_1}
-y_{21}X^{-s_1\omega_1} 
+y_{42}X^{-\omega_2},
\cr
[s_1s_2s_1s_2] &=
(1-y_{32})
+(y_{22}+y_{42}+y_{43}+y_{53})X^{-\omega_1}
-y_{32}X^{-s_1\omega_1}
-y_{32}X^{-s_2s_1\omega_1}
\cr &\qquad\qquad\qquad\qquad
-y_{32}X^{-\omega_2}
-y_{32}X^{-s_2\omega_2},
\cr
[s_2s_1s_2] &=
(1-y_{21}+y_{42})
+(y_{42}-y_{21}-y_{52}-y_{53}-y_{63})X^{-\omega_1}
+(y_{42}-y_{21})X^{-s_1\omega_1}
\cr &\qquad\qquad\qquad
+(y_{42}-y_{21})X^{-s_2s_1\omega_1}
+y_{42}X^{-\omega_2}
+y_{42}X^{-s_2\omega_2},
\cr
[s_1s_2s_1] &=
(1-2y_{32})
+(y_{22}+y_{42}+y_{43}+y_{53})X^{-\omega_1}
+(y_{22}+y_{42}+y_{43}+y_{53})X^{-s_1\omega_1}
\cr &\qquad\qquad
-y_{32}X^{-s_2s_1\omega_1}
-y_{32}X^{-s_1s_2s_1\omega_1}
\cr &\qquad\qquad
-(y_{32}+y_{43}+y_{53})X^{-\omega_2}
-y_{32}X^{-s_2\omega_2}
-y_{32}X^{-s_1s_2\omega_2},
\cr
[s_2s_1] &=
(1-y_{21}+2y_{42})
+(y_{42}-y_{21}-y_{52}-y_{53}-y_{63})X^{-\omega_1}
\cr &\qquad
+(y_{42}-y_{21}-y_{32}-y_{53}-y_{63})X^{-s_1\omega_1}
+(y_{42}-y_{21})X^{-s_2s_1\omega_1}
\cr &\qquad
+(y_{42}-y_{21})X^{-s_1s_2s_1\omega_1}
+(y_{42}+y_{63})X^{-\omega_2}
+y_{42}X^{-s_2\omega_2}
+y_{42}X^{-s_1s_2\omega_2},
\cr
[s_1s_2] &=
1-y_{11}-y_{21}-y_{32}-y_{43}-y_{53}
+(y_{22}+y_{32})(1+y_{10}+y_{20})X^{-\omega_1}
\cr &\qquad
+(y_{22}+y_{32}+y_{42})X^{-s_1\omega_1}
+(y_{22}+y_{32}+y_{42})X^{-s_2s_1\omega_1}
\cr &\qquad
-(y_{32}+y_{43}+y_{53})X^{-\omega_2}
-(y_{32}+y_{43}+y_{53})X^{-s_2\omega_2}
-y_{32}X^{-s_1s_2\omega_2}
-y_{32}X^{-s_2s_1s_2\omega_2},
\cr
[s_2] &=
(1+y_{31}+y_{32}+2y_{42}+y_{63})
-(y_{21}+y_{52}+y_{53}+y_{84})X^{-\omega_1}
-(y_{21}+y_{52}+y_{53})X^{-s_1\omega_1}
\cr &\qquad
-(y_{21}+y_{52}+y_{53})X^{-s_2s_1\omega_1}
-y_{21}X^{-s_1s_2s_1\omega_1}
-y_{21}X^{-s_2s_1s_2s_1\omega_1}
\cr &\qquad
+(y_{42}+y_{63})X^{-\omega_2}
+(y_{42}+y_{63})X^{-s_2\omega_2}
+y_{42}X^{-s_1s_2\omega_2}
+y_{42}X^{-s_2s_1s_2\omega_2},
\cr
[s_1] &=
1-(y_{11}+y_{21}+y_{32}+2y_{43}+2y_{53})
+(y_{22}+y_{54})(1+y_{10}+y_{20})X^{-\omega_1}
\cr &\qquad
+(y_{22}+y_{54})(1+y_{10}+y_{20})X^{-s_1\omega_1}
+(y_{22}+y_{32}+y_{42})X^{-s_2s_1\omega_1}
\cr &\qquad
+(y_{22}+y_{32}+y_{42})X^{-s_1s_2s_1\omega_1}
-(y_{32}+y_{43}+y_{53}+y_{64})X^{-\omega_2}
-(y_{32}+y_{43}+y_{53})X^{-s_2\omega_2}
\cr &\qquad
-(y_{32}+y_{43}+y_{53})X^{-s_1s_2\omega_2}
-y_{32}X^{-s_2s_1s_2\omega_2}
-y_{32}X^{-s_1s_2s_1s_2\omega_2},
\cr
[1] &=
(1+y_{31}+y_{42}+y_{63}-y_{53}-y_{43})
-y_{21}(1+y_{32})^2X^{-\omega_1}
\cr &\qquad
+y_{22}(1+y_{10}+y_{20})(1+y_{21}+y_{31})X^{-s_1\omega_1}
-(y_{21}+y_{52}+y_{53})X^{-s_2s_1\omega_1}
\cr &\qquad
+y_{22}X^{-s_1s_2s_1\omega_1}
-y_{21}X^{-s_2s_1s_2s_1\omega_1}
-y_{32}(1+y_{11})(1+y_{21})X^{-\omega_2}
+(y_{42}+y_{63})X^{-s_2\omega_2}
\cr &\qquad
-(y_{32}+y_{43}+y_{53})X^{-s_1s_2\omega_2}
+y_{42}X^{-s_2s_1s_2\omega_2}
-y_{32}X^{-s_1s_2s_1s_2\omega_2}
+y_{53}X^{-\rho}.
}$$
The multiplication of the Schubert classes is given by
%
%Multiplication by $[1]$:
%
$$\eqalign{
[1]^2 
&=
\alpha_{10}\alpha_{01}\alpha_{11}
\alpha_{21}\alpha_{31}\alpha_{32}[1],
\cr
[1][s_1]
&=
-\alpha_{01}\alpha_{11}\alpha_{21}
\alpha_{31}\alpha_{32}[1],
\cr
[1][s_2]
&=
-\alpha_{10}\alpha_{11}\alpha_{21}
\alpha_{31}\alpha_{32}[1],
\cr
[1][s_1s_2]
&=
\alpha_{11}\alpha_{21}\alpha_{31}\alpha_{32}[1],
\cr
[1][s_2s_1]
&=
\alpha_{11}\alpha_{21}\alpha_{31}\alpha_{32}[1],
\cr
[1][s_1s_2s_1]
&=
-\alpha_{11}\alpha_{21}\alpha_{32}(1+y_{11}+y_{21})[1], 
}
\qquad\qquad
\eqalign{
[1][s_2s_1s_2]
&=
-\alpha_{21}\alpha_{31}\alpha_{32}[1], \cr
[1][s_1s_2s_1s_2]
&=
\alpha_{21}\alpha_{32}(1+y_{21})[1], \cr
[1][s_2s_1s_2s_1]
&=
\alpha_{21}\alpha_{32}(1+y_{21})[1], \cr
[1][s_1s_2s_1s_2s_1]
&=
-\alpha_{32}(1+y_{32})[1], \cr
[1][s_2s_1s_2s_1s_2]
&=
-\alpha_{21}(1+y_{21})[1], \cr
\cr
}$$
%
%\bigskip
%%Multiplication by $[s_1]$:
%
$$\eqalign{
[s_1]^2
&=
-\alpha_{01}\alpha_{11}\alpha_{21}\alpha_{31}\alpha_{32}[s_1]
\cr
[s_1][s_2]
&=
\alpha_{11}\alpha_{21}\alpha_{31}\alpha_{32}[1]
\cr
[s_1][s_1s_2]
&=
-\alpha_{11}\alpha_{21}\alpha_{32}(y_{01}+y_{11}+y_{21})[1]
+\alpha_{01}\alpha_{11}\alpha_{21}\alpha_{32}[s_1]
\cr
[s_1][s_2s_1]
&=
\alpha_{11}\alpha_{21}\alpha_{31}\alpha_{32}[s_1]
\cr
[s_1][s_1s_2s_1]
&=
-\alpha_{11}\alpha_{21}\alpha_{32}(1+y_{11}+y_{21})[s_1] 
\cr
[s_1][s_2s_1s_2]
&=
\alpha_{21}\alpha_{32}(y_{11}+y_{21})[1]
-\alpha_{11}\alpha_{21}\alpha_{32}[s_1] 
\cr
[s_1][s_1s_2s_1s_2]
&=
-\alpha_{32}(y_{22}+y_{32})[1]
+\alpha_{11}\alpha_{32}(1+y_{11})[s_1]
\cr
[s_1][s_2s_1s_2s_1]
&=
\alpha_{21}\alpha_{32}(1+y_{21})[s_1]
\cr
[s_1][s_1s_2s_1s_2s_1]
&=
-\alpha_{32}(1+y_{32})[s_1]
\cr
[s_1][s_2s_1s_2s_1s_2]
&=
y_{32}[1]-\alpha_{32}[s_1]
\cr
}$$
%
%\bigskip
%%Multiplication by $[s_2]$:
%
$$\eqalign{
[s_2]^2
&=
-\alpha_{10}\alpha_{11}\alpha_{21}\alpha_{31}\alpha_{32}[s_2] 
\cr
[s_2][s_1s_2]
&=
\alpha_{11}\alpha_{21}\alpha_{31}\alpha_{32}[s_2] 
\cr
[s_2][s_2s_1]
&=
-\alpha_{21}\alpha_{31}\alpha_{32}y_{10}[1]
+\alpha_{10}\alpha_{21}\alpha_{31}\alpha_{32}[s_2]
\cr
[s_2][s_1s_2s_1]
&=
\alpha_{21}\alpha_{32}(y_{21}+y_{31})[1]
-\alpha_{21}\alpha_{31}\alpha_{32}[s_2] 
\cr
[s_2][s_2s_1s_2]
&=
-\alpha_{21}\alpha_{31}\alpha_{32}[s_2]
\cr
[s_2][s_1s_2s_1s_2]
&=
\alpha_{21}\alpha_{32}(1+y_{21})[s_2]
\cr
[s_2][s_2s_1s_2s_1]
&=
-\alpha_{21}(y_{31}+y_{52})[1]
+\alpha_{21}\alpha_{31}(1+y_{21})[s_2]
\cr
[s_2][s_1s_2s_1s_2s_1]
&=
y_{63}[1]
-\alpha_{21}(1+y_{21}+y_{42})[s_2]
\cr
[s_2][s_2s_1s_2s_1s_2]
&=
-\alpha_{21}(1+y_{21})[s_2]
}$$
%
%\bigskip
%%Multiplication by $[s_1s_2]$:
%
$$\eqalign{
[s_1s_2]^2
&=
-\alpha_{11}\alpha_{21}\alpha_{32}(y_{01}+y_{11}+y_{21})[s_2]
+\alpha_{01}\alpha_{11}\alpha_{21}\alpha_{32}[s_1s_2]
\cr
[s_1s_2][s_2s_1]
&=
\alpha_{21}\alpha_{32}(y_{11}+y_{21}+\alpha_{31})[1]
-\alpha_{11}\alpha_{21}\alpha_{32}[s_1]
-\alpha_{21}\alpha_{31}\alpha_{32}[s_2]
\cr
[s_1s_2][s_1s_2s_1]
&=
-\alpha_{32}( y_{32}+y_{42}
\{\,+\alpha_{11}(y_{01}+2y_{11}+y_{21})\,\})
[1]
+\alpha_{11}\alpha_{32}(y_{01}+y_{11})[s_1]
\cr &\qquad
+\big(
\alpha_{31}\alpha_{32}y_{11}
+\alpha_{11}\alpha_{32}(y_{01}+y_{11}+y_{21})
\big)[s_2]
-\alpha_{01}\alpha_{11}\alpha_{32}[s_1s_2] 
\cr
[s_1s_2][s_2s_1s_2]
&=
\alpha_{21}\alpha_{32}(y_{11}+y_{21})[s_2]
-\alpha_{11}\alpha_{21}\alpha_{32}[s_1s_2]
\cr
[s_1s_2][s_1s_2s_1s_2]
&=
-\alpha_{32}(y_{22}+y_{32})[s_2]
+\alpha_{11}\alpha_{32}(1+y_{11})[s_1s_2]
\cr
[s_1s_2][s_2s_1s_2s_1]
&=
\big(
y_{63}\,\{+\alpha_{32}(y_{11}+y_{21})\,\}
\big)[1]
-\alpha_{32}y_{11}[s_1]
-\big(
\alpha_{32}(y_{11}+y_{21})
+\alpha_{31}y_{32}
\big)[s_2]
\cr &\qquad
+\alpha_{11}\alpha_{32}[s_1s_2]
\cr
[s_1s_2][s_1s_2s_1s_2s_1]
&=
\{\,-(y_{33}+y_{43}+y_{53})[1]\,\}
+y_{33}[s_1]
+(y_{33}+y_{43}+y_{53})[s_2]
\cr &\qquad
-\alpha_{11}(1+y_{11}+y_{22})[s_1s_2]
\cr
[s_1s_2][s_2s_1s_2s_1s_2]
&=
y_{32}[s_2]-\alpha_{32}[s_1s_2]
}$$
%
%\bigskip
%%Multiplication by $[s_2s_1]$:
%
$$\eqalign{
[s_2s_1]^2
&=
-\alpha_{21}\alpha_{31}\alpha_{32}y_{10}[s_1]
+\alpha_{10}\alpha_{21}\alpha_{31}\alpha_{32}[s_2s_1]
\cr
[s_2s_1][s_1s_2s_1]
&=
\alpha_{21}\alpha_{31}(y_{21}+y_{31})[s_1]
-\alpha_{21}\alpha_{31}\alpha_{32}[s_2s_1]
\cr
[s_2s_1][s_2s_1s_2]
&=
-\alpha_{21}(y_{51}+y_{52}
\{\,+\alpha_{31}y_{10}\,\})[1]
+\alpha_{21}(\alpha_{10}y_{31}+\alpha_{32}y_{10})[s_1]
\cr &\qquad
+\alpha_{21}\alpha_{31}(y_{10}+y_{21})[s_2]
-\alpha_{10}\alpha_{21}\alpha_{31}[s_2s_1]
\cr
[s_2s_1][s_1s_2s_1s_2]
&=
\big(
y_{62}\{\,+\alpha_{31}(y_{21}+y_{31})\,\}
\big)[1]
-\big(
\alpha_{31}y_{21} +\alpha_{10}(y_{31}+y_{41})
\big)[s_1]
\cr &\qquad
-\big(
\alpha_{31}y_{21}+\alpha_{32}y_{31}
\big)[s_2]
+\alpha_{21}\alpha_{31}[s_2s_1]
\cr
[s_2s_1][s_2s_1s_2s_1]
&=
-\alpha_{21}(y_{31}+y_{52})[s_1]
+\alpha_{21}\alpha_{31}(1+y_{21})[s_2s_1]
\cr
[s_2s_1][s_1s_2s_1s_2s_1]
&=
y_{63}[s_1]
-\alpha_{21}(1+y_{21}+y_{42})[s_2s_1]
\cr
[s_2s_1][s_2s_1s_2s_1s_2]
&=
\{\,-y_{31}[1]\,\}
+y_{31}[s_1]
+y_{31}[s_2]
-\alpha_{31}[s_2s_1]
}$$

\bigskip
%Multiplication by $[s_1s_2s_1]$:

$$\eqalign{
[s_1s_2s_1]^2
&=
-\alpha_{32}(y_{32}+y_{42}
\{\,+\alpha_{11}(y_{11}+y_{21})\,\}) [s_1]
\cr &\qquad
+\big(\alpha_{11}\alpha_{32}(y_{01}+y_{11}+y_{21})
+\alpha_{31}\alpha_{32}y_{11}\big)[s_2s_1]
-\alpha_{01}\alpha_{11}\alpha_{32}[s_1s_2s_1] 
\cr
[s_1s_2s_1][s_2s_1s_2]
&=
\big(1\{\,+\alpha_{11}(y_{11}+y_{22}+y_{33}+y_{31}+y_{42})
+\alpha_{31}(y_{21}+y_{32})
+\alpha_{32}y_{21}\,\}\big) [1] 
\cr &\qquad
-\big(\alpha_{11}(y_{21}+\alpha_{32})
+\alpha_{10}(y_{31}+y_{41}+y_{32}+y_{42})\big)[s_1]
\cr &\qquad
-(\alpha_{31}(y_{21}+y_{32})+\alpha_{11}(y_{21}+y_{32}+y_{31}+\alpha_{42})[s_2]
\cr &\qquad
+\alpha_{11}\alpha_{32}[s_1s_2]
+\alpha_{21}\alpha_{31}[s_2s_1]
\cr
[s_1s_2s_1][s_1s_2s_1s_2]
&=
\{\,-(y_{33}+2y_{43}+y_{53}
+\alpha_{11}(y_{01}+y_{11})
+\alpha_{21}(y_{11}+y_{21})
)[1]\,\}
\cr &\qquad
+\big( 
y_{33}+y_{43}
\{\,+\alpha_{11}(y_{01}+y_{11})
+\alpha_{21}(y_{11}+y_{21})\,\}
\big) [s_1]
\cr &\qquad
\big(
(y_{33}+y_{43}+y_{53})
\{\,+\alpha_{11}(y_{01}+y_{11})
+\alpha_{21}(y_{11}+y_{21})\,\}
\big)[s_2]
\cr &\qquad
-\alpha_{11}(y_{01}+y_{11}+y_{22})
[s_1s_2]
-\big(
\alpha_{11}(y_{01}+y_{11})
+\alpha_{21}(y_{11}+y_{21})
\big)[s_2s_1]
\cr &\qquad
+\alpha_{01}\alpha_{11} [s_1s_2s_1]
\cr
[s_1s_2s_1][s_2s_1s_2s_1]
&=
(y_{62}\{\,+\alpha_{32}y_{21}\,\})[s_1]
-\big(\alpha_{31}y_{32}+\alpha_{32}(y_{11}+y_{21})\big)[s_2s_1]
+\alpha_{11}\alpha_{32} [s_1s_2s_1]
\cr
[s_1s_2s_1][s_1s_2s_1s_2s_1]
&=
\{\,-(y_{43}+y_{53})[s_1]\,\}
+(y_{33}+y_{43}+y_{53})[s_2s_1]
-\alpha_{11}(1+y_{11}+y_{22}) [s_1s_2s_1]
\cr
[s_1s_2s_1][s_2s_1s_2s_1s_2]
&=
\{\,(y_{11}+y_{21})[1]
-(y_{11}+y_{21})[s_1]
-(y_{11}+y_{21})[s_2] \,\}
\cr &\qquad
+y_{11}[s_1s_2]
+(y_{11}+y_{21})[s_2s_1]
-\alpha_{11}[s_1s_2s_1]
\cr
}
$$
%
%\bigskip
%%Multiplication by $[s_2s_1s_2]$:
%
$$\eqalign{
[s_2s_1s_2]^2
&=
-\alpha_{21}(y_{21}+y_{42})[s_2]
+\big(
\alpha_{11}\alpha_{21}y_{31} +\alpha_{21}\alpha_{31}y_{10}
\big)[s_1s_2]
-\alpha_{10}\alpha_{21}\alpha_{31}[s_2s_1s_2]
\cr
[s_2s_1s_2][s_1s_2s_1s_2]
&=
y_{53}[s_2]
-\big(
\alpha_{21}y_{31}
+\alpha_{11}\alpha_{21}\alpha_{32}y_{21}
\big)[s_1s_2]
+\alpha_{21}\alpha_{31}[s_2s_1s_2] 
\cr
[s_2s_1s_2][s_2s_1s_2s_1]
&=
\{\,-\big(
y_{51}+y_{52}
+\alpha_{31}y_{10}
\big)[1]\,\}
+(y_{41}\{\,+\alpha_{31}y_{10}\,\})[s_1]
+(y_{42}+y_{52}\{\,+\alpha_{31}y_{10}\,\})[s_2]
\cr &\qquad
-(\alpha_{11}y_{31}+\alpha_{31}y_{10})[s_1s_2]
-\alpha_{31}y_{10} [s_2s_1]
+\alpha_{10}\alpha_{31}[s_2s_1s_2]
\cr
[s_2s_1s_2][s_1s_2s_1s_2s_1]
&=
\{\,(y_{31}+y_{32}+y_{42})[1]
-(y_{31}+y_{32})[s_1]
-(y_{31}+y_{32}+y_{42})[s_2]\,\}
\cr &\qquad
+(y_{31}+y_{32})[s_1s_2]
+y_{31}[s_2s_1]
-\alpha_{31}[s_2s_1s_2]
\cr
[s_2s_1s_2][s_2s_1s_2s_1s_2]
&=
y_{31}[s_1s_2]
-\alpha_{31}[s_2s_1s_2]
}$$
%
%\bigskip
%%Multiplication by $[s_1s_2s_1s_2]$:
%
$$\eqalign{
[s_1s_2s_1s_2]^2
&=
\{\,-y_{43}[s_2]\,\}
+(y_{32}+y_{42}\{\, +\alpha_{01}y_{21}+\alpha_{32}y_{11}\,\})[s_1s_2]
\cr &\qquad
-\big(\alpha_{01}(y_{11}+y_{21})+\alpha_{31}(y_{01}+y_{11})
\big)[s_2s_1s_2]
+\alpha_{01}\alpha_{11}[s_1s_2s_1s_2]
\cr
[s_1s_2s_1s_2][s_2s_1s_2s_1]
&=
\{\,(y_{21}+y_{31}+y_{32}+y_{42}+\alpha_{11})[1]
\cr &\qquad
-(y_{21}+y_{31}+y_{32}+\alpha_{11})[s_1] 
-(y_{21}+y_{31}+y_{32}+y_{42}+\alpha_{11})[s_2]\,\}
\cr &\qquad
+(y_{31}+y_{42}\{,+\alpha_{11}\,\})[s_1s_2]
+(y_{21}+y_{31}\{\,+\alpha_{11}\,\})[s_2s_1]
\cr &\qquad
-\alpha_{11}[s_1s_2s_1]
-\alpha_{31}[s_2s_1s_2]
\cr
[s_1s_2s_1s_2][s_1s_2s_1s_2s_1]
&=
\{\,-(y_{01}+y_{11}+y_{21}+y_{22}+y_{32})[1]
\cr &\qquad
+(y_{01}+y_{11}+y_{21}+y_{22})[s_1]
+(y_{01}+y_{11}+y_{21}+y_{22}+y_{32})[s_2]
\cr &\qquad
-(y_{01}+y_{11}+y_{21}+y_{22})[s_1s_2]
-(y_{01}+y_{11}+y_{21})[s_2s_1]\,\}
\cr &\qquad
+y_{01}[s_1s_2s_1]
+(y_{01}+y_{11}+y_{21})[s_2s_1s_2]
-\alpha_{01}[s_1s_2s_1s_2]
\cr
[s_1s_2s_1s_2][s_2s_1s_2s_1s_2]
&=
\{\,-y_{21}[s_1s_2]\,\}
+(y_{11}+y_{21})[s_2s_1s_2]
-\alpha_{11}[s_1s_2s_1s_2]
}$$
%
%\bigskip
%%Multiplication by $[s_2s_1s_2s_1]$:
%
$$\eqalign{
[s_2s_1s_2s_1]^2
&=
\{\,-y_{52}[s_1]
+(y_{42}+y_{52})[s_2s_1]\,\}
-(\alpha_{11}y_{31}+\alpha_{31}y_{10})[s_1s_2s_1]
+\alpha_{10}\alpha_{31}[s_2s_1s_2s_1]
\cr
[s_2s_1s_2s_1][s_1s_2s_1s_2s_1]
&=
\{\,y_{42}[s_1]
-(y_{31}+y_{41})[s_2s_1]\,\}
+(y_{31}+y_{32})[s_1s_2s_1]
-\alpha_{31}[s_2s_1s_2s_1]
\cr
[s_2s_1s_2s_1][s_2s_1s_2s_1s_2]
&=
\{\,-y_{10}[1]
+y_{10}[s_1]
+y_{10}[s_2]
-y_{10}[s_1s_2]
-y_{10}[s_2s_1]\,\}
\cr &\qquad
+y_{10}[s_1s_2s_1]
+y_{10}[s_2s_1s_2]
-\alpha_{10}[s_2s_1s_2s_1]
}$$
%
%\bigskip
%%Multiplication by $[s_1s_2s_1s_2s_1]$:
%
$$\eqalign{
[s_1s_2s_1s_2s_1]^2
&=
\{\,-y_{32}[s_1]
+(y_{22}+y_{32})[s_2s_1]
-(y_{11}+y_{21}+y_{22})[s_1s_2s_1]\,\}
\cr &\qquad
+(y_{01}+y_{11}+y_{21})[s_2s_1s_2s_1]
-\alpha_{01}[s_1s_2s_1s_2s_1]
\cr
[s_1s_2s_1s_2s_1][s_2s_1s_2s_1s_2]
&=
\{\,[1]-[s_1]-[s_2]+[s_1s_2]+[s_2s_1]
\cr &\qquad
-[s_1s_2s_1]-[s_2s_1s_2]\,\}
+[s_1s_2s_1s_2]+[s_2s_1s_2s_1]
}$$
%
%\bigskip
%%Multiplication by $[s_2s_1s_2s_1s_2]$:
%
$$\eqalign{
[s_2s_1s_2s_1s_2]^2
&=
y_{10}[s_1s_2s_1s_2]
-\alpha_{10}[s_2s_1s_2s_1s_2]
}$$

%\vfill\eject
\bigskip\bigskip

\section 5. References

\bigskip

\medskip
\item{[BGG]} {\smallcaps I.N.\ Bernstein, I.M.\ Gel'fand and S.I.\ Gel'fand},
{\it Schubert cell and cohomology of the spaces $G/P$}, 
Russ.\ Math.\ Surv.\ {\bf 28} (3) (1973), 1--26.

\medskip
\item{[Br]} {\smallcaps M.\ Brion},
{\it Positivity in the Grothendieck group of complex flag varieties},
J.\ Alg.\ {\bf 258} no.\ 1 (2002), 137--159.

\medskip
\item{[Bou]} {\smallcaps N.\ Bourbaki},
{\sl Groupes et algebres de Lie}, Chapt.\ IV-VI,
Masson, Paris, 1981.

\medskip
\item{[Ch]} {\smallcaps C.\ Chevalley}, 
{\it Sur les decompositions cellulaires des espaces $G/B$},  
in {\sl Algebraic Groups and their Generalizations:
Classical Methods}, W.\ Haboush and B.\ Parshall eds.,
Proc.\ Symp.\ Pure Math., Vol.\ {\bf 56} Pt.\ 1, 
Amer.\ Math.\ Soc.\ (1994), 1--23.

\medskip
\item{[CG]} {\smallcaps N.\ Chriss and V.\ Ginzburg}, 
{\sl Representation theory and complex geometry}, Birkh\"auser, 
Boston, 1997.

\medskip
\item{[D]} {\smallcaps M.\ Demazure}, 
{\it D\'esingularisation des vari\'et\'es de Schubert g\'en\'eralis\'ees}, 
Ann.\ Sci.\ \'Ecole Norm.\ Sup.\ {\bf 7} (1974), 53--88.

\medskip
\item{[FL]} {\smallcaps W.\ Fulton and A.\ Lascoux},
{\it A Pieri formula in the Grothendieck ring of a flag bundle},
Duke Math.\ J.\ {\bf 76} (1994), 711--729.

\medskip
\item{[Fu]} {\smallcaps W.\ Fulton},
{\sl Intersection Theory},  Ergebnisse der Mathematik (3) {\bf 2},
Springer-Verlag, Berlin-New York, 1984. 

\medskip
\item{[Gd]} {\smallcaps A.\ Grothendieck},
{\it Sur quelques propri\'et\'es fondamnetales en th\'eorie des
intersections}, in {\sl Anneaux de Chow et applications},
S\'eminaire C.\ Chevalley 2e ann\'ee
(mimeographed notes) Paris (1958), pages 4-01~--~4-36.

\medskip
\item{[Gr]} {\smallcaps W.\ Graham},
{\it Positivity in equivariant Schubert calculus}, Duke Math.\ J.\ 
{\bf 109} (2001), 599--614.

\medskip
\item{[Hz]} {\smallcaps F.\ Hirzebruch},
{\sl Topological methods in algebraic geometry},
Third edition, Springer-Verlag, 1995.

\medskip
\item{[KK]} {\smallcaps B.\ Kostant and S.\ Kumar}, 
{\it $T$-equivariant K-theory of generalized flag varieties}, 
J.\ Differential Geom.\ {\bf 32} (1990), 549--603. 

\medskip
\item{[KR]} {\smallcaps C.\ Kriloff and A.\ Ram}, 
{\it Representations of graded Hecke algebras},
Representation Theory {\bf 6} (2002), 31--69.

\medskip
\item{[La]} {\smallcaps A.\ Lascoux}, 
{\it Chern and Yang through ice},
preprint 2002.

\medskip
\item{[L1]} {\smallcaps P.\ Littelmann},
{\it A Littlewood-Richardson rule for symmetrizable Kac-Moody algebras},
Invent.\ Math.\ {\bf 116} (1994), 329-346.

\medskip
\item{[L2]} {\smallcaps P.\ Littelmann},
{\it Paths and root operators in representation theory}, Ann.\ Math.\ 
{\bf 142} (1995), 499-525.

\medskip
\item{[L3]} {\smallcaps P.\ Littelmann},
{\it Characters of representations and paths in ${\goth H}_{\RR}^*\;$},
Proc.\ Symp.\ Pure Math.\ {\bf 61} (1997), 29-49.

\medskip
\item{[LS]} {\smallcaps P.\ Littelmann and C.S.\ Seshadri},
{\it A Pieri-Chevalley formula for $K(G/B)$ and standard monomial
theory}, in {\sl Studies in memory of Issai Schur}, Progress in Mathematics 210, Birkh\"auser, 2003, 155-176.

%\medskip
%\item{[Mac]} {\smallcaps I.G.\ Macdonald},
%{\sl Algebraic geometry: Introduction to schemes}, W.\ A.\ Benjamin, New
%York-Amsterdam,1968.

\medskip
\item{[Ma]} {\smallcaps O.\ Mathieu},
{\it Positivity of some intersections in $K_0(G/B)$},
in {\sl Commutative algebra, homological algebra and representation 
theory (Catania/Genoa/Rome, 1998)}, J.\ Pure Appl.\ Algebra 
{\bf 152} (2000), no.\ 1-3, 231--243. 

\medskip
\item{[NR]} {\smallcaps K.\ Nelsen and A.\ Ram},
{\it Kostka-Foulkes polynomials and Macdonald spherical functions},
in {\sl Surveys in Combinatorics 2003}, C.\ Wensley ed., London Math.\ 
Soc.\  Lect. Notes {\bf 307} 
Camb. Univ.\ Press (2003), 325--370.

\medskip
\item{[P]} {\smallcaps H.\ Pittie}, 
{\it Homogeneous vector bundles over homogeneous spaces},
Topology {\bf 11} (1972), 199--203.

\medskip
\item{[PR1]} {\smallcaps H.\ Pittie and A.\ Ram},
{\it A Pieri-Chevalley formula in the K-theory of a $G/B$ bundle},
Elec.\ Research Announcements {\bf 5} (1999), 102--107.

\medskip
\item{[PR2]} {\smallcaps H.\ Pittie and A.\ Ram},
{\it A Pieri-Chevalley formula in the K-theory of flag variety},
preprint 1998, 
{\tt http://www.math.wisc.edu/\~{}ram/preprints.html}. 

\medskip
\item{[R]} {\smallcaps A.\ Ram},
{\it Affine Hecke algebras and generalized standard Young tableaux},
J.\ Algebra {\bf 260} (2003), 367--415.

%\medskip
%\item{[Sn]} {\smallcaps V.\ Snaith},
%{\it ???}, ???.

\medskip
\item{[St]} {\smallcaps R.\ Steinberg},
{\it On a theorem of Pittie}, Topology {\bf 14} (1975), 173--177.

\vfill\eject
\end

\lemma Let $w_0$ be the longest element of $W$ and let
$w_i$ be the longest element of the (parabolic) subgroup 
$W_i={\rm Stab}(\omega_i)$.  Let $w^i\in W$ be such that $w_0 = w_iw^i$ and 
$\ell(w_0)=\ell(w_i)+\ell(w^i)$.
Then
$$T_{w_i}\varepsilon_{w^i} = (-1)^{\ell(w^i)}(T_{w_0}-T_{s_iw_0}).$$
\pf
Since $T_{w_0}$ is, up to constant multiples, the unique
element of $H$ such that $T_jT_{w_0} = T_{w_0}$ for all $1\le j\le n$
it is sufficient to establish that 
$$T_j(T_{s_iw_0}+(-1)^{\ell(w^i)}T_{w_i}\varepsilon_{w^i}) = 
T_{s_iw_0}+(-1)^{\ell(w^i)}T_{w_i}\varepsilon_{w^i},
\qquad\hbox{for $1\le j\le n$.}
$$
If $j\ne i$ then $s_j\in W_i$ and 
$$T_j(T_{s_iw_0}+(-1)^{\ell(w^i)}T_{w_i}\varepsilon_{w^i})
=T_jT_{s_iw_0}+(-1)^{\ell(w^i)}T_jT_{w_i}\varepsilon_{w^i}
=T_{s_iw_0}+(-1)^{\ell(w_i)}T_{w_i}\varepsilon_{w^i},$$
since $s_jw_i < w_i$ and the unique $k$ such that
$s_k(s_iw_0)> s_iw_0$ is $k=i$.  Since ?????
$$T_i(T_{s_iw_0}+(-1)^{\ell(w^i)}T_{w_i}\varepsilon_{w^i})
=T_{w_0}+(-1)^{\ell(w^i)}T_iT_{w_i}\varepsilon_{w^i}
=T_{s_iw_0}+(T_{w_0}-T_{s_iw_0}+(-1)^{\ell(w^i)}T_iT_{w_i}\varepsilon_{w^i})
=?????
$$
which completes the proof.
\endpf